\documentclass[AMA,STIX1COL]{WileyNJD-v2}

\usepackage{amssymb,amsmath}
\usepackage{bbm}%,mathabx}
\usepackage{amsfonts}
\usepackage{mathrsfs,mathtools}

\makeatletter
\def\BState{\State\hskip-\ALG@thistlm}
\makeatother

\newcommand{\bs}[1]{\boldsymbol{#1}}
\newcommand{\bmu}{\bs{\mu}}

\newcommand{\bx}{\bs{x}}

\newcommand{\R}{\mathbbm{R}}

\newcommand{\calD}{\mathcal{D}}
\newcommand{\calJ}{\mathcal{J}}

\newcommand{\calN}{\mathcal{N}}

\DeclareMathOperator*{\argmax}{argmax}

\articletype{Research article}%

\received{}
\revised{}
\accepted{}

\begin{document}

\title{Adaptive greedy algorithms based on parameter-domain decomposition and reconstruction for the reduced basis method}
%\protect\thanks{Sth.}}

\author[1]{Jiahua Jiang}

\author[2]{Yanlai Chen*}

\authormark{Jiang, Chen: \textsc{Multi-level greedy algorithms for RBM}}

\address[1]{\orgdiv{Department of Mathematics}, \orgname{Virginia Tech University}, \orgaddress{\state{VA}, \country{U.S.A.}}\\ Email: {\tt{jiahua@vt.edu}}}

\address[2]{\orgdiv{Department of Mathematics}, \orgname{University of Massachusetts Dartmouth}, \orgaddress{\state{MA}, \country{U.S.A.}} \email{\tt{Yanlai.Chen@umassd.edu}}}

\corres{*Corresponding author.}

\abstract[Summary]{The reduced basis method (RBM) empowers repeated and rapid evaluation of parametrized partial differential equations through an offline-online decomposition, a.k.a. a learning-execution process. A key feature of the method is a greedy algorithm repeatedly scanning the training set, a fine discretization of the parameter domain, to identify the next dimension of the parameter-induced solution manifold along which we expand the surrogate solution space. Although successfully applied to problems with fairly high parametric dimensions, the challenge is that this scanning cost dominates the offline cost due to it being proportional to the cardinality of the training set which is exponential with respect to the parameter dimension. In this work, we review three recent attempts in effectively delaying this curse of dimensionality, and propose two new hybrid strategies through successive refinement and multilevel maximization of the error estimate over the training set. All five offline-enhanced methods and the original greedy algorithm are tested and compared on {two types of problems: the thermal block problem and the geometrically parameterized Helmholtz problem.}
}


\maketitle

\section{Introduction}

The reduced basis method (RBM) \cite{Quarteroni2015, HesthavenRozzaStammBook} has proved to be a viable option for the purpose of designing fast numerical algorithms  for parametrized systems.  The parameters {(denoted by $\bmu$ throughout this paper)} may include boundary conditions, material properties, amount of uncertainty, geometric configurations, source properties etc. To describe a realistic system, practitioners often resort to a large number of parameters making the the generation of a reduced model {computationally challenging}.

The {critical ingredient enabling RBM to attain} orders of magnitude gain in marginal (i.e. per parameter instance) computation time, is an offline-online decomposition process where the basis selection (i.e. surrogate solution space building) is performed offline by a greedy algorithm, see e.g. \cite{Rozza2008, Haasdonk2017Review, Prudhomme_Rovas_Veroy_Maday_Patera_Turinici, NguyenVeroyPatera2005} for details. The ultimate goal is that the complexity of the reduced solver, presumably to be called for an overwhelming number of times or even in a realtime fashion, is independent of the degrees of freedom {(denoted by $\calN$ throughout this paper)} of the {full order model (FOM, a.k.a. ``truth'' in the RB literature)} approximation.

{The construction of} the {reduced basis} space 
relies on a greedy scheme that keeps track of an  efficiently-computable error estimator/indicator, denoted by $\Delta_n(\bmu)$, indicating the discrepancy between the dimension-$n$ surrogate solution {(RB solution)} and the {FOM} solution. This greedy procedure starts by selecting the first parameter $\bmu^1$ randomly from the training set {(i.e. a discretized version of the parameter domain)} and obtaining
its corresponding truth approximation $u^\mathcal{N}(\bmu^1)$. 
We then have a one-dimensional RB space $W_1 = \{u^{\mathcal N}(\bmu^1)\}$. Next, the scheme 
obtains an RB approximation ${u}^\calN_{1}(\boldsymbol{\mu})$ for each parameter in the training set together with $\Delta_1(\bmu)$. The greedy choice for the $(n+1)$th parameter $(n=1,2, \cdots)$ is made and the RB space augmented by
\begin{equation}
\label{eq:rbmgreedy}
\bmu^{n+1} = \underset{{\bmu} \in \Xi_{\rm train}}{\argmax} \Delta_{n}({\bmu}), \quad W_{n+1} = W_n \oplus \{u^\calN(\bmu^{n+1})\}.
\end{equation}

The other parts of the offline process are devoted to the necessary preparations for the online reduced solver, a variational (i.e. Galerkin or Petrov-Galerkin) projection \cite{benner2015survey, CarlbergBouMoslehFarhat2011, CarlbergBaroneAntil2017} into the surrogate space, or a {representation} of the solution at strategically chosen points in the physical domain \cite{ChenGottlieb2013, CGJMX2019}. The greedy procedure \eqref{eq:rbmgreedy} implies that the offline cost is proportional to  the cardinality of the training set $\Xi_{\rm train}$. It is essential for this training set to be fine enough so that we are not missing critical phenomena not represented by $\Xi_{\rm train}$. In particular, its cardinality is exponential with respect to the parameter dimension. 
As a consequence, the ``break-even'' number of simulations for the parametric system (i.e. minimum number of simulations that makes the offline preparation stage worthwhile) is exponential with respect to the parameter dimension as well, severely diminishing the attractiveness of RBM.

The need of delaying {this} curse of dimensionality motivates numerous recent attempts in designing offline-enhancement strategies for systems with high-dimensional parameter domain. Indeed, the authors of \cite{sen2008reduced} propose {a static decomposition of} the {given parameter training set $\Xi_{\rm train}$} {\it a priori} into a sequence of subsets that increase geometrically in size, and then perform the classical greedy algorithm on each of them {sequentially. The idea is that scanning smaller sets will quickly result in a reduced solver capable of resolving a significant part of the larger sets making them more affordable for scanning. }Authors of \cite{hesthaven2014efficient} take a different route. They fix the cardinality of an active training set {that is} much smaller than the full training set. {After each greedy step, the active training set is pruned (with the parameters whose error estimates are below the tolerance) and replenished to the original cardinality with random selections from the full training set.} Finally, \cite{jiang2017offline} promotes to recycle the information afforded by the error estimates $\{\Delta_n(\bmu)\}$ and to {adaptively} construct a surrogate training set. {The adaptive nature of the construction which, moreover, is informed by the error estimate makes it potentially more effective than the approach of \cite{sen2008reduced}.} 
There have been other techniques in the literature to alleviate the RB offline cost such as the parameter domain adaptivity \cite{haasdonk2011training, iapichino2016reduced}, greedy sampling acceleration through nonlinear optimization \cite{urban2014greedy, iapichino2015optimization}, and local reduced basis method \cite{ohlberger2015error, kaulmann2015localized, ohlberger2017true,zou2019adaptive}, etc.

In this paper, we focus on the three approaches of \cite{sen2008reduced, hesthaven2014efficient, jiang2017offline} that are similar in nature. We review the essence of each algorithm, {compare them,} and more importantly, propose two new hybrid strategies through successive refinement and multilevel maximization of $\{\Delta_n(\bmu)\}$ over the training set. By testing all these five strategies and the classical greedy, we demonstrate that the two new hybrids outperform others some of which may even fail for {a} more challenging {geometrically parameterized} Helmholtz {problem}.

The paper is organized as follows. In Section \ref{sec:background}, we review the {three} existing offline enhancement approaches {mentioned above}. {In Section \ref{sec:multi-greedy}, the two new hybrid approaches are introduced.} Numerical results for two {test} problems to demonstrate the accuracy and efficiency of these offline {improvement} methods are shown in Section \ref{sec: numerics}. Finally, concluding remarks are drawn in Section \ref{sec5}.

\section{Background}
\label{sec:background}
In this section,  we briefly review the classical reduced basis method and the three offline enhancement strategies. This also serves as a background and motivation for the subsequent discussion of our newly proposed multi-level greedy techniques leading to more efficient reduced basis methods. Table \ref{tab:notation} outlines our notations. 

\begin{table}
  \begin{center}
  \resizebox{\textwidth}{!}{
    \renewcommand{\tabcolsep}{0.4cm}
    \renewcommand{\arraystretch}{1.3}
    {\scriptsize
    \begin{tabular}{@{}lp{0.8\textwidth}@{}}
      \toprule
      $\bmu$ & Parameter in $\calD \subseteq \R^p$ \\
      $u(\bmu)$ & Function-valued solution of a parameterized PDE \\
      $\mathcal{N}$ & Degrees of freedom (DoF) in PDE ``truth" solver \\
      $u^{\mathcal{N}}(\bmu)$ & Truth solution (finite-dimensional)\\
      $N$ & Number of reduced basis snapshots, $N \ll \mathcal{N}$\\
      $N_{max}$ & Terminal number of reduced bases \\
      $\bmu^j$ & ``Snapshot" parameter values, $j=1, \ldots, N$\\
      $X^{\mathcal{N}}_{{N}}$ & Span of $u^{\mathcal{N}}\left(\bmu^k\right)$ for $k=1, \ldots, N$\\
      $u_N^{\mathcal{N}}(\bmu)$ & Reduced basis solution, $u_N^{\mathcal{N}} \in X^{\mathcal{N}}_{{N}}$\\
      $e_N(\bmu)$ & Reduced basis solution error, equals $u^{\mathcal{N}}(\bmu) - u_N^{\mathcal{N}}(\bmu)$ \\
      $\Xi_{\rm{train}}$ & Parameter training set, a finite subset of $\mathcal{D}$ \\
      $n_{\rm train}$ & Size of $\Xi_{\rm{train}}$\\
      $\Delta_{{N}} \left(\bmu\right)$ & Error estimate (upper bound) for $\left\|e_N\left(\bmu\right)\right\|$ \\
      $\epsilon_{\mathrm{tol}}$ & Error estimate stopping tolerance in greedy sweep \\
      \midrule
      & {\textbf{A priori training set decomposition}} \\
      $\mathcal{J}$ & Number of sample sets\\
      $\Xi_{\rm train, j}$ & The $j-$th sample set\\
      $\Delta_{N,j}(\bmu)$ & A posteriori error estimate with $N$ reduced bases at $j-$th sample set\\
      \midrule
      & {\textbf{Adaptive enriching}} \\
      $\Xi$ & Sample parameters of fixed cardinality $(<< n_{\rm train})$ \\
      {$M_{\rm sample}$} & Size of the active training set $\Xi$, much smaller than $n_{\rm train}$\\
      \midrule
      & {\textbf{Adaptive construction of surrogate training sets}} \\
      $\ell$ & Number of ``outer" loops in the offline enhancement procedure, Algorithm \ref{alg:oerbm} \\
      $E_{\ell}$ & The largest error estimator at the beginning of outer loop $\ell$\\
      $\Xi_{\mathrm{sur}}$ & Surrogate Training Set (STS), a subset of $\Xi_{\rm{train}}$  \\
      $K_{\mathrm{damp}}$ & A constant integer controlling the damping ratio of the error estimate when examining the STS\\ 
      {$C_M$} & {A constant integer adjusting the size of surrogate training set}\\
    \bottomrule
    \end{tabular}
  }
    \renewcommand{\arraystretch}{1}
    \renewcommand{\tabcolsep}{12pt}
  }
  \end{center}
  \caption{Notation used throughout this article.}\label{tab:notation}
\end{table}

\subsection{Reduced Basis Methods}

{We first provide an overview of the essential ingredients of RBM in its classical form. For more details see e.g. \cite{Quarteroni2015, HesthavenRozzaStammBook,Rozza2008, Haasdonk2017Review}. Indeed,} given ${\bmu} \in \mathcal{D}$, the goal is to evaluate a certain output of interest
\begin{equation}
\label{eq:interest}
s(\bmu) = \ell(u(\bmu); \bmu){,}
\end{equation}
where the function $u(\bmu) \in X$ satisfies
\begin{equation}\label{eq:sat}
a(u(\bmu),v; \bmu) = f(v; \bmu), \quad v \in X,
\end{equation}
which is a parametrized partial differential equation (pPDE) written in a weak form with $\bmu \in \calD$ being the (possibly multi-dimension) parameter. Here $X = X(\Omega)$ is a Hilbert space satisfying, e.g., $H^{1}_{0}(\Omega) \subset X(\Omega) \subset H^{1}(\Omega)$. We denote by $(\cdot , \cdot)_{X}$ the inner product associated with the space $X$, whose induced norm $|| \cdot ||_{X} = \sqrt{(\cdot , \cdot)_{X}}$ is equivalent to the usual $H^{1}(\Omega)$ norm. 
We {typically} assume that $a(\cdot ,\cdot; \bmu): X \times X \rightarrow \mathbb{R}$ is continuous and coercive over $X$ uniformly in $\bmu \in \mathcal{D}$, that is,
\begin{subequations}
\begin{equation}\label{eq:con}
\gamma(\bmu) \coloneqq \underset{w \in X}{\sup} ~ \underset{v \in X}{\sup}~ \frac{a(w, v; \bmu)}{||w||_{X}||v||_{X}} < \infty, \quad \forall \bmu \in \mathcal{D}, 
\end{equation} 
\begin{equation}\label{eq:coe}
\alpha(\bmu) \coloneqq \underset{w \in X} \inf \frac{|a(w, w; \bmu)|}{||w||^2_{X}} \geqslant \alpha_{0} > 0, \forall \bmu \in \mathcal{D}. 
\end{equation}
\end{subequations}
$f(\cdot)$ and $\ell(\cdot)$ are linear continuous functionals over $X$. 
We assume that there is a finite-dimensional discretization for the model problem \eqref{eq:sat}: The solution space $X$ is discretized by an $\calN$-dimensional space $X^{\mathcal{N}}$ (i.e., $\dim(X^{\mathcal{N}}) = \mathcal{N}$) and \eqref{eq:interest} and \eqref{eq:sat} are discretized as
\begin{equation}\label{eq:update_problem}
\begin{cases}
\text{For } \bmu \in \mathcal{D},~ \text{solve} \\
s^{\mathcal{N}} = \ell(u^{\mathcal{N}}(\bmu); \bmu) ~ \text{where} ~ u^{\mathcal{N}}(\bmu) \in X^{\mathcal{N}} ~ \text{satisfies} \\
a(u^{\mathcal{N}}, v; \bmu) = f(v; \bmu) \quad \forall v \in X^{\mathcal{N}}.
\end{cases}
\end{equation}
The relevant quantities such as the coercivity constant \eqref{eq:coe} are defined according to the discretization,
\[
\alpha^{\mathcal{N}}(\bmu)  = \inf_{w \in X^{\mathcal{N}}} \frac{a(w, w; \bmu)}{||w||^{2}_{X}} \ge \alpha(\bmu) \ge \alpha_0, ~ \forall \bmu \in \mathcal{D}.
\]
In the RBM literature, any quantity associated to $\mathcal{N}$ is called a ``truth" {(FOM)}. E.g., $u^{\mathcal{N}}$ is called the ``truth solution", \eqref{eq:update_problem} ``truth solver''. $\calN$ is typically very large so that resolving the FOM gives highly accurate approximations for all $\bmu \in \mathcal{D}$.

For a training set $\Xi_{\rm train} \subset \mathcal{D}$ which consists of a fine discretization of $\mathcal{D}$ of finite cardinality and a collection of $N$ parameters $S_{N} = \{\bmu^{1}, \dots, \bmu^{N}\} \subset \Xi_{\rm train}$, we define the reduced basis space as 
\begin{equation}\label{eq:space_define}
X^{\mathcal{N}}_{N} \coloneqq \text{span} \{ u^{\mathcal{N}}(\bmu^{n}), 1 \leq n \leq N \}.
\end{equation}
The reduced basis approximation is now defined as: Given $\bmu \in \mathcal{D}$, seek a surrogate RB solution $u_{N}^{\mathcal{N}}(\bmu)$ by solving the following {reduced} system
\begin{equation}
\label{eq:reduced_system}
\begin{cases}
\text{For } \bmu \in \mathcal{D},~ \text{evaluate} \\
s_{N}^{\mathcal{N}} = \ell(u_{N}^{\mathcal{N}}(\bmu); \bmu) ~ \text{s.t.} ~ u_{N}^{\mathcal{N}}(\bmu) \in X_{N}^{\mathcal{N}} \subset X^{\mathcal{N}}~ \text{satisfies} \\
a(u_{N}^{\mathcal{N}}, v; \bmu) = f(v) \quad \forall v \in X_{N}^{\mathcal{N}}.
\end{cases}
\end{equation}
\eqref{eq:reduced_system} is called the reduced solver, i.e. Reduced Order Model (ROM). 
The typical multiple orders of magnitude speedup of RBM manifests from that the assembly of ROM is independent of $\calN$, which crucially relies on the affine assumption of the parameter dependent problem \eqref{eq:sat}, such as, 
\begin{equation}
\label{eq:affine}
a(w, v; \bmu) = \sum_{q = 1}^{Q_{a}} \Theta_{a}^{q}(\bmu)a^{q}(w,v), \quad \text{and} \quad f(v; \bmu) = \sum_{q = 1}^{Q_{f}} \Theta_{f}^{q}(\bmu)f^{q}(v).
\end{equation} 
Here $\Theta_{a}^{q}$, $\Theta_{f}^{q}$ {are} $\bmu-$dependent functions, and $a^{q}$, $f^{q}$ are $\bmu-$independent forms. {With Galerkin procedure in (\ref{eq:reduced_system}), the RB solution $u_{N}^{\mathcal{N}}$ for any parameter $\bmu \in \mathcal{D}$ can be represented as 
\begin{equation}
\label{eq:rb-sol}
u_N^{\mathcal{N}}(\bmu) = \sum_{i=1}^{N}u^{\mathcal{N}}_{Ni}(\bmu)u^{\mathcal{N}}(\bmu^i),
\end{equation}
where $\{u^{\mathcal{N}}_{Ni}(\bmu)\}_{i=1}^{N}$ are the RB coefficients obtained by solving (\ref{eq:reduced_system}) \footnote{{In actual computation, we orthonormalize the basis $\{u^{\mathcal{N}}(\bmu^i)\}$ through a Gram-Schmidt procedure for numerical stability.}}.
}
With {the affine} hypothesis {(\ref{eq:affine})}, we can apply an offline-online decomposition to enable fast resolution of the ROM \eqref{eq:reduced_system}. {In the} offline stage, {we replace the reduced basis solution in (\eqref{eq:reduced_system}) by (\ref{eq:rb-sol}) and choose $v = u^{\mathcal{N}}(\bmu^j), 1\le j \le N$ as our test functions. Then we obtain the RB ``stiffness'' equations
\begin{equation}
\label{eq:stiff}
\sum_{i=1}^{N}a(u^\calN(\bmu^i), u^\calN(\bmu^j);\bmu)u^{\mathcal{N}}_{Ni}(\bmu) = f(u^\calN(\bmu^j);\bmu) \quad 1 \le j \le N.
\end{equation}
Toward a fast construction of the stiff matrix, we plug \eqref{eq:rb-sol} into \eqref{eq:stiff} to obtain,
\begin{equation}
\sum_{q=1}^{Q_a}\sum_{i=1}^{N}\Theta_{a}^{q}(\bmu)a^{q}(u^\calN(\bmu^i), u^\calN(\bmu^j))u^{\mathcal{N}}_{Ni}(\bmu) = \sum_{q=1}^{Q_f}\Theta_{f}^{q}(\bmu)f^{q}(u^\calN(\bmu^j)) \quad 1 \le j \le N.
\end{equation}
and precompute} $a^{q}(u^\calN(\bmu^i), u^\calN(\bmu^j))$ and {$f^{q}(u^\calN(\bmu^j))$}, which {are} relatively expensive but only done once. In the online stage, we construct the matrices and vectors in the reduced system (\ref{eq:reduced_system}) and solve the resulting reduced basis problem. We remark that when assumption {(\ref{eq:affine})} is violated, we turn to Empirical Interpolation Method (EIM) \cite{Barrault_Nguyen_Maday_Patera, Grepl_Maday_Nguyen_Patera, ChaturantabutSorensen2010} to approximate the non-affine operators by affine ones. This is the case for our second test problem {in Section \ref{sec:helm}}.

Leading to the key error estimator $\Delta_N(\bmu)$, we define the error $e_N(\bmu) := u^{\mathcal{N}}(\bmu) - u^{\mathcal{N}}_{N}(\bmu) \in X^{\mathcal{N}}$. {The} linearity of $a(\cdot, \cdot; \bmu)$ yields the following error equation:
\begin{equation}
\label{eq:residual_eq}
a(e_N(\bmu), v; \bmu) = r_N(v; \bmu) \quad \forall v \in X^\calN,
\end{equation}
where the residual $r_N(\cdot; \bmu) \in (X^{\mathcal{N}})'$ (the dual of $X^{\calN}$) operated on $v \in X^\calN$ is defined as $f(v; \bmu) - a(u_{N}^{\mathcal{N}}(\bmu), v; \bmu)$.  We define the {\it a posteriori} error estimator for the solution, {which is a rigorous bound for the error,} as
\begin{equation}\label{eq:error_estimator_en}
  \Delta_{N}(\bmu) = \frac{\lVert r_N(\cdot; \bmu)\rVert_{(X^{\mathcal{N}})'}}{\alpha^{\mathcal{N}}_{LB}(\bmu)} {\ge \lVert e_N(\bmu) \rVert_{(X^{\mathcal{N}})},}
\end{equation}
where $\alpha^{\mathcal{N}}_{LB}(\bmu)$ is a lower bound of the coercivity constant {$\alpha^{\mathcal{N}}(\bmu)$. The otherwise expensive evaluation of $\alpha^{\mathcal{N}}(\bmu)$ for all $\bmu$ can be done efficiently using the so-called Successive Constraint Method  and related approaches \cite{huynh2007successive, huynh2010natural, chen2016certified, yano2014space, manzoni2015heuristic}. As to the dual norm of the residual $r_N(\cdot; \bmu) \in (X^{\mathcal{N}})'$, we also take advantage of a suitable offline-online splitting. In fact, we first invoke the Riesz representation theorem and define functions $\mathcal{C}^{q}, \mathcal{L}_i^{q} \in X^{\mathcal{N}}$ such that 
\begin{equation}\label{error_es_problem}
\begin{cases}
  (\mathcal{C}^{q}, v)_{X^\calN} = f^{q}(v)_{X^{\calN}} \quad \forall v \in X^{\mathcal{N}}, 1 \le q \le Q_f \\
  (\mathcal{L}_{i}^{q}, v)_{X^{\calN}} = a^{q}(u^{\calN}\left(\bmu^i\right), v) \quad \forall v \in X^{\mathcal{N}}, 1\le q \le Q_a, 1 \le i \le N.
\end{cases}
\end{equation}
Notice that problem (\ref{error_es_problem}) is parameter-independent, hence $\mathcal{C}^{q}$ and $\mathcal{L}_{m}^{q}$ can be computed offline. Combining (\ref{eq:residual_eq}), (\ref{eq:rb-sol}) and (\ref{eq:affine}), we have
\begin{align}\nonumber
  \lVert r_N(\cdot; \bmu)\rVert_{(X^{\mathcal{N}})'}^2 &= 
\sum_{{q_1} = 1}^{Q_f}  \sum_{{\tilde{q}_1} = 1}^{Q_f} \Theta_f^{q_1}(\bmu) \Theta_f^{{\tilde{q}_1}}(\bmu) (\mathcal{C}^{{q_1}}, \mathcal{C}^{{\tilde{q}_1}})_{X^{\calN}} + \\
&\sum_{q_2 = 1}^{Q_{a}}\sum_{i = 1}^{N}\Theta_a^{q_2}(\bmu)u^{\mathcal{N}}_{Ni} \left\{ \sum_{q_2' = 1}^{Q_{a}}\sum_{i' = 1}^{N}\Theta_a^{q_2'}(\bmu)u^{\mathcal{N}}_{Ni'}(\mathcal{L}_{i}^{q_2}, \mathcal{L}_{i'}^{q_2'})_{X^{\calN}}\right\} \nonumber\\
\label{eq:error_es_quan}
&- 2\sum_{q_2 = 1}^{Q_{a}}\sum_{i = 1}^{N} {\sum_{q_1 = 1 }^{Q_f}} \Theta_a^{q_2}(\bmu)u^{\mathcal{N}}_{Ni}(\bmu)({\mathcal{C}^{q_1}}, \mathcal{L}_{i}^{q_2})_{X^{\calN}},
\end{align}
where the parameter-independent quantities $(\mathcal{C}^{q_1}, \mathcal{C}^{{\tilde{q}_1}})_{X^{\calN}}, (\mathcal{C}{^{q_1}},\mathcal{L}_{i}^{q_2})_{X^{\calN}}, (\mathcal{L}_{i}^{q_2}, \mathcal{L}_{i'}^{q_2'})_{X^{\calN}}, \,\, 1 \le i, i' \le N_{\rm RB}, 1 \le {{q_1}}, {{\tilde{q}_1}} \le Q_{f}, 1 \le q_2, q_2' \le Q_{a}$ can be precomputed and stored. 
Once we can efficiently calculate $\Delta_n(\bmu)$, } the classical greedy algorithm, outlined in Algorithm \ref{alg:subgreedy} is invoked to build the parameter set $S_{N}$ and the {resulting} reduced basis space $X^{\mathcal{N}}_{N}$.
\begin{algorithm}
\begin{algorithmic}[1]
\If{$N = 0$}
\State Initialization: Choose an initial parameter value $\bmu^1 \in \Xi_{\rm train}$, set $S_1 = \{\bmu^1\}$, compute $u^\calN(\bmu^{1})$, and {let $X^\calN_N = \{u^\calN(\bmu^1)\}$} $N=1$;
\EndIf
\While{$\underset{\bmu \in {\Xi_{\rm train}}}{\max}\Delta_{N}(\bmu) > \varepsilon_{\mathrm{tol}}$}
\State Choose $\bmu^{N+1} =  \underset{\bmu \in {\Xi_{\rm train}}}{\argmax} \Delta_{N}(\bmu)$;
\State $S_{N+1} = S_{N}\bigcup \{\bmu^{N+1}\}$;
\State Compute $u^\calN(\bmu^{N+1})$ and augment RB space $X^\calN_{N+1} = X^\calN_N \oplus \{u^\calN(\bmu^{N+1})\}$.
\State $N \gets N+1$.
\EndWhile
\end{algorithmic}
 \caption{Classical Greedy, $\left(N, X_N^\calN\right) = {\mathbb {CG}}(\Xi_{\rm train}, \varepsilon_{\mathrm{tol}}, N, X_N^\calN)$} 
 \label{alg:subgreedy} 
\end{algorithm}

\subsection{Existing offline enhancement strategies}
{The} greedy algorithm ${\mathbb {CG}}(\Xi_{\rm train}, \varepsilon_{\mathrm{tol}}, N, X_N^\calN)$ requires maximization of the {\it a posteriori} error estimate over $\Xi_{\rm train}$. This becomes a bottleneck in the construction of the RB space $X_N^\calN$, especially when the parameter domain $\calD$ is of high dimension. Much recent research has {developed} {the} schemes capable of accelerating this {standard} greedy algorithm. In this section, we review three such {offline enhancement} strategies as proposed in \cite{sen2008reduced, hesthaven2014efficient, jiang2017offline} in Subsection \ref{sec:senRB} to \ref{sec:smm-greedy} respectively. The improved greedy algorithms are described in Algorithms \ref{sen-greedy} to \ref{alg:oerbm} accordingly.

\subsubsection{A priori training set decomposition (TSD)}
\label{sec:senRB}
In \cite{sen2008reduced}, a modified greedy algorithm is provided to address the many-parameter heat conduction problems. 
It attempts to {mitigate the effect of large $\Xi_{\rm train}$ by running} the classical greedy algorithm {first on a relatively small training set before attempting to sequentially augment the reduced basis space in larger training sets.} {To do so, they decompose the full training set $\Xi_{\rm train}$ into} a sequence of mutually exclusive subsets \{$\Xi_{\rm train, 1}, \dots, \Xi_{\rm train, \calJ}\}$ {that has  \{$\Xi_{\rm train, 1}, \dots, \Xi_{\rm train, \calJ - 1}\}$ increasing in size geometrically and the last one ensuring that the whole training set is covered, i.e. $\Xi_{\rm train, 1} \cup \dots \cup \Xi_{\rm train, \calJ} = \Xi_{\rm train}$}.
The {main} idea is { running the greedy algorithm over the small sets would result in a basis capable of resolving a large portion of the larger sets. These portions, i.e.} the samples whose current error estimator is below the tolerance, {will be skipped thereby increasing the speed of the greedy scan. The TSD algorithm includes a final run on $\Xi_{\rm train}$ as a sanity check.}

The size of the {first} sample set $n_{\rm tr, small} = |\Xi_{\rm train, 1}|$ serves as the {(only)} tuning parameter for this approach. As an example to build up this {partition}, let  $\mathcal{J} = {\rm floor}(\log_2 (\frac{n_{\rm train}}{ n_{\rm tr, small}})$. For $j = 2: \mathcal{J}-1$, $|\Xi_{\rm train, j}| = 2^{j-1}n_{\rm tr, small}$. The size of $\Xi_{\rm train, \mathcal{J}}$ equals to $ n_{\rm train} - \sum_{j=1}^{\mathcal{J}-1}|\Xi_{\rm train, j}|$. To make each of the subsets span the whole training set, the algorithm randomly sample $n_{\rm train, j}$ points from $\Xi_{\rm train}$. For $j$-th sample set, the {\it a posteriori} error estimate computed by $N$ reduced bases is denoted as $\Delta_{N, j}(\bmu)$. The pseudo-code of this methodology is provided in algorithm \ref{sen-greedy}.

\begin{remark}
This approach has $n_{\rm tr, small}$ as the only tuning parameter.
Small $n_{\rm tr, small}$ likely leads to higher dimension of the surrogate {solution} space negatively impacting the online efficiency. 
Therefore, we usually choose relatively large $n_{\rm tr, small}$ in comparison to $n_{\rm train}$. Unfortunately, large $n_{\rm tr, small}$ will lead to more costly offline stage. This is obvious since 
when, {in the extreme case}, $n_{\rm tr, small} = n_{\rm train}$, Algorithm \ref{sen-greedy} is exactly the same as classical greedy algorithm. {This balance of tuning the parameter $n_{\text{tr, small}}$ is rather intricate, as shown by our second numerical experiment in Section \ref{sec:diff}. In comparison, our proposed adaptive approaches are much less sensitive to their tuning parameters.}
\end{remark}

\begin{algorithm}
\begin{algorithmic}[1]
\State {Determine a partition $\{\Xi_{\rm train, j}\}_{j=1}^\calJ$ for $\Xi_{\rm train}$, e.g. by randomly sampling $n_{\rm train, j}$ points from $\Xi_{\rm train}$ to form $\Xi_{\rm train, j}$ where $|\Xi_{\rm train, j}| = 2^{j-1}n_{\rm tr, small}$ for $j = 1, \dots, \mathcal{J}-1$, and  $\Xi_{\rm train, \calJ}$ completing the partition.}
\State {Let $N = 0$.}
\For{$j = 1: \mathcal{J}$}
\State Call $\left(N, X_N^\calN\right) = {\mathbb {CG}}(\Xi_{\rm train, j}, \varepsilon_{\mathrm{tol}}, N, X_N^\calN)$
\EndFor
\end{algorithmic}
 \caption{Training set decomposition based classical greedy ${\left(N, X_N^\calN\right) = \mathbb {TSD\_CG}}(\Xi_{\rm train}\, \varepsilon_{\mathrm{tol}})$}
 \label{sen-greedy}
\end{algorithm}

\subsubsection{Adaptive enriching}
\label{sec: AE}
{To avoid running standard greedy algorithm multiple times on the large training set, the adaptive enriching algorithm \cite{hesthaven2014efficient} opts for executing the greedy algorithm on a dynamically determined subset of the full training set that is fixed in size much smaller than the full training set. }
{The sample} set is {iteratively} updated, after each greedy step, by removing parameters that have error estimate below the tolerance and {randomly} adding new parameter values from the training set. {The rationale is that} it is not worthwhile to keep those parameters whose corresponding surrogate solutions are already accurate enough.  {The size of the sample set is always maintained at a fixed number $M_{\rm sample}$ that needs to be specified by practitioners. } 
More details of this algorithm are provided in  Algorithm \ref{Jan-greedy}. Since the maximum {\it a posteriori} error estimate in the stopping criteria is only the maximum in {sample set rather than full training set}, {a safety check step} is introduced to {ensure} each parameter of the full training set is checked once. For problems that we don't have monotonic decay of error or error estimate for, a final ``check'' over the entire training set is necessary.

\begin{remark}
The adaptive enriching scheme has one tuning parameter $M_{\rm sample}$. Two extremal cases are worth mentioning:
\begin{itemize}
    \item When $M_{\rm sample} = 1$, the method becomes the approach taken by \cite{elman2013reduced} and \cite{LiuChenChenShu2019} where each parameter is examined once and decision on keep or toss made immediately. This approach may result in a larger-than-necessary surrogate space.
    \item When $M_{\rm sample} = n_{\rm train}$, this algorithm is identical to the classical greedy algorithm, thus no savings are realized.
\end{itemize}
{Furthermore, a tuning parameter test of this algorithm for the many-parameter heat conduction problem is presented in Section \ref{sec:diff}, which indicates that the performance of this approach is again rather sensitve to the choice of its tuning parameter $M_{\rm sample}$.}

\end{remark}

\begin{algorithm}
\begin{algorithmic}[1]
\State $N_{safe} = {\rm ceil}(|\Xi_{\rm train}|/M_{\rm sample})$;
\State Randomly generate an initial training set $\Xi$ with $M_{\rm sample}$ parameters from $\Xi_{\rm train}$;
\State Choose an initial parameter $\bmu^{1} \in \Xi$.  Set $S_{1} = \{ \bmu^{1}\}$, $X^\calN_{1} = {\rm span} \{u^\calN(\bmu^{1}) \}$ and $N = 1$;
\State Set ${\rm safe} =0$, $\varepsilon = 2\varepsilon_{\mathrm{tol}}$ and $r = n_{\rm train}$;
\While{($\varepsilon > \varepsilon_{\mathrm{tol}} $ or ${\rm safe} \le N_{\rm safe}$) and ($r > 0$)}
\State Choose $\bmu^{N+1} =  \underset{\bmu \in \Xi}{\argmax} \Delta_{N}(\bmu)$;
\State Augment RB space $X^\calN_{N+1} = X^\calN_N \oplus \{u^\calN(\bmu^{N+1})\}$, $S_{N+1} = S_{N}\bigcup \{\bmu^{N+1}\}$;
\State Truncate $\Xi$ by $\Xi_{<\varepsilon} = \{ \bmu \in \Xi: \Delta_N(\bmu) < \varepsilon_{\mathrm{tol}}\}$.
\State Truncate $\Xi_{\rm train}$ by $\Xi_{<\varepsilon}$, and set $r = r -  |\Xi_{<\varepsilon}|$;
\If{$|\Xi_{<\varepsilon}| = M_{\rm sample}$}
\State Set ${\rm safe} =  {\rm safe} +1$.
\EndIf
\State Randomly choose $M_{\rm sample} - |\Xi|$ parameters from $\Xi_{\rm train}$ for addition to $\Xi$;
\State Set $N \gets N+1$;
\EndWhile
\end{algorithmic}
 \caption{Adaptive enriching classical greedy $\left(N, X_N^\calN\right) = {\mathbb {AE\_CG}}(\Xi_{\rm train}, \varepsilon_{\mathrm{tol}}, M_{\rm sample})$}
  \label{Jan-greedy}
\end{algorithm}

\subsubsection{Adaptive construction of surrogate training sets}
\label{sec:smm-greedy}
{The offline-enhanced RBMs of  \cite{jiang2017offline} adaptively} identifies a subset of the training set, {termed} a ``Surrogate Training Set'' (STS), {on which to} perform {the classical} greedy algorithm. {Its distinctive feature is that the construction of STS is informed by the error estimator $\{\Delta_n(\bmu): \bmu \in \Xi_{\rm train}\}$ {while the methods in \cite{sen2008reduced, hesthaven2014efficient} are not}. After every sweep of the full parameter domain $\Xi_{\rm train}$ {via performing $\bmu^{N+1} =  \underset{\bmu \in {\Xi_{\rm train}}}{\argmax} \Delta_{N}(\bmu)$}, {a much smaller STS $\Xi_{\rm sur}$ is constructed} and then the classical greedy is invoked on the STS until it is deemed as fully resolved, {i.e. $\underset{\bmu \in {\Xi_{\rm sur}}}{\max} \Delta_{N}(\bmu)$ is below the tolerance}. At that point, $\Xi_{\rm train}$ is examined again {by performing the maximization over $\Xi_{\rm train}$} to start the next cycle.  These surrogate sets are significantly smaller than the full training set, yet selected well enough to capture the general landscape of the {error-estimate} manifold. Thus, one attraction of this approach is its unique algorithmic structure, shown in Algorithm \ref{alg:oerbm}, of error estimate-informed transition between the full and surrogate training sets. 
The frequent targeted searches over STS produce a potential relative saving of $1 - |\Xi_{\rm sur}|/|\Xi_{\rm train}|$.  The reason is that these searches are operated on STS $\Xi_{\rm sur}$ instead of on $\Xi_{\rm train}$, {and that the cost of these searches is linearly dependent on the size of the set it operates on}.}

{
We adopt the more cost effective of the two strategies from \cite{jiang2017offline} for building STS, the 
Successive Maximization Method (SMM). SMM is inspired by inverse transform sampling for non-standard univariate probability distributions. 
}
Indeed, motivated by the notion that the difference between the norm of the errors $\lVert e(\bmu_1) \rVert_{X} - \lVert e(\bmu_2) \rVert_{X} $ is partially indicative of the difference between the solutions. When selecting the {$(N+1)$}-th parameter value, {we construct STS by equidistantly sampling values from $\Delta_N(\Xi_{\rm train})$.} 
{With} $\epsilon_{\rm tol}$ the stopping tolerance for the RB sweep, let $\Delta_N^{\rm max} = \underset{{\bmu \in \Xi_{\rm train}}}{\max} \Delta_N(\bmu)$. We define $I_N^{M_{\ell}}$ as an {equi-spaced} set between $\epsilon_{\rm tol}$ and $\Delta_N^{\rm max}$:
\begin{align} \label{eq:partition}
  I_N^{M_{\ell}} = \left\{ {\nu_{N,m} \coloneqq \,} \epsilon_{\rm tol} + (\Delta_N^{\rm max} - \epsilon_{\rm tol}) \frac{m}{M_{\ell}}: m = {0}, \dots, {M_{\ell} - 1}\right\}.
\end{align}
Roughly speaking, we attempt to construct $\Xi_{\rm Sur}$ as $\Xi_{\rm Sur} = \Delta_N^{-1}\left(I_N^{M_{\ell}}\right) \bigcap \Xi_{\rm train}$. 
On outer loop round $\ell$, we have $|\Xi_{\rm Sur}| {{\le}} M_{\ell}$ by this construction, where $M_{\ell}$ can be chosen as any monotonically increasing function with respect to $\ell$. But in order to avoid excessively large $\Xi_{sur}$, we set $M_{\ell} = C_{M}(\ell + 1)$, where $C_{M}$ is a constant. After this construction, we repeatedly sweep the current  $\Xi_{\rm Sur}$ until 
    \begin{equation*}
      \max_{\bmu \in \Xi_{\rm Sur}} \Delta_N(\bmu) \le E_{\ell}\, \frac{1}{((\ell+1) \times K_{\rm damp})},
    \end{equation*}
    where $E_{\ell}$ is the starting (global) maximum error estimate for this outer loop iteration. The damping ratio for outer loop $\ell$, {$\frac{1}{(\ell+1) \times K_{\rm damp}}$}, enforces that the maximum error estimate over the $\Xi_{\rm Sur}$ decreases by a damping factor controlled by $K_{\rm damp}$ which should be determined by the practitioner and the problem at hand.  Following the choice in \cite{jiang2017offline}, we take $K_{\rm damp}$ to be constant in this paper. Algorithm \ref{alg:oerbm} details the  adaptively constructed STS approach.

\begin{remark} 
This approach has two tuning parameters $C_{M}$ and $K_{\rm damp}$. The choice of $C_M$ indirectly controls the size of $\Xi_{\rm sur}$. These surrogate training sets should be small enough compared to the full one $\Xi_{\rm train}$ to offer considerable acceleration of the greedy sweep, yet large enough to capture the general landscape of the solution manifold that is iteratively learned. On the other hand, $K_{damp}$ controls how {accurate} we intend to resolve {$u^{\mathcal{N}}_{N}(\bmu)$, where $\bmu \in \Xi_{\rm sur}$}.  The larger $C_{M}$ and $K_{\rm damp}$ are, the faster algorithm \ref{alg:oerbm} will be. 
{
As will be shown by the tuning parameter study in Section \ref{sec:diff}, the efficiency of this approach is less sensitive to the choice of its tuning parameters. 
}
\end{remark}

\begin{algorithm}
\begin{algorithmic}[1]
\If{$N=0$}
\State Randomly select the first sample $\bmu^1 \in \Xi_{\rm train}$, and set $ {N} = 1$, $\ell = 0$ and $E_0 = 2 \varepsilon_{\mathrm{tol}}$;
\State Obtain truth solution $u^\mathcal{N}(\bmu^1)$, and set $S_{1} = \{ \bmu^{1}\}$, $X^\calN_1 = \mbox{span}\left\{u^{\mathcal N}(\bmu^1)\right\}$;
\EndIf
\While {$(E_{\ell} > \varepsilon_{\mathrm{tol}})$}
\vspace{0.05in}
\State Set $\ell \leftarrow \ell + 1$;
\State 
\begin{minipage}{0.03\textwidth}
\rotatebox{90}{One-step greedy}
\end{minipage} 
\begin{minipage}{0.03\textwidth}
\rotatebox{90}{scan on $\Xi_{\rm train}$}
\end{minipage} 
\hspace{4mm}
%one-step greedy starts
\framebox[0.8 \textwidth]
{\begin{minipage}{0.8 \textwidth}
\For{each $\boldsymbol{\mu} \in \Xi_{\rm train}$}
\State Obtain RBM solution $u^{\mathcal{N}}_{ {N}}(\boldsymbol{\mu}) \in X^\calN_{ {N}}$ and error estimate ${\Delta_{ {N}}}(\boldsymbol{\mu})$;
\EndFor
\vspace{0.05in}
\State $\bmu^{ {N}+1} = \underset{\bmu \in \Xi_{\rm train}}{\argmax} \Delta_{ {N}}(\bmu)$, $\varepsilon = \Delta_{ {N}}(\bmu^{ {N}+1})$, $E_{\ell} = \varepsilon$;
\State Augment RB space $X^\calN_{ {N}+1} = X^\calN_{ {N}} \oplus \{u^\calN(\bmu^{ {N}+1})\}$;
\State $S_{ {N}+1} = S_{ {N}}\bigcup \{\bmu^{ {N}+1}\}$;
\State Set $ {N} \leftarrow  {N}+1$;
\end{minipage}
}
%one-step greedy ends
\vspace{0.05in}
\State \hspace{-0.0\textwidth} Construct STS $\Xi_{\rm Sur}$ based on $\{(u_{ {N}-1}^\calN(\bmu), \Delta_{ {N}-1}(\bmu)): \bmu \in \Xi_{\rm train}\}$ with SMM;
\vspace{0.05in}
\State 
\begin{minipage}{0.03\textwidth}
\rotatebox{90}{Multi-step greedy}
\end{minipage} 
\begin{minipage}{0.03\textwidth}
\rotatebox{90}{scan on $\Xi_{\rm Sur}$}
\end{minipage} 
\hspace{4mm}
%multi-step greedy starts
\framebox[0.8 \textwidth]
{\begin{minipage}{0.8 \textwidth}
\While {$(\varepsilon > \varepsilon_{\mathrm{tol}})$ and {$(\frac{\varepsilon}{E_{\ell}} >  \, \frac{1}{K_{\rm damp} \times (\ell + 1)})$}}
\vspace{0.05in}
\For{each $\boldsymbol{\mu} \in \Xi_{\rm Sur}$}
\State Obtain RBM solution $u^{\mathcal{N}}_{ {N}}(\boldsymbol{\mu}) \in X^\calN_{ {N}}$ and error estimate ${\Delta_{ {N}}}(\boldsymbol{\mu})$;
\EndFor
\vspace{0.05in}
\State $\bmu^{ {N}+1} = \underset{\bmu \in \Xi_{\rm Sur}}{\argmax} \Delta_{ {N}}(\bmu)$, $\varepsilon = \Delta_{ {N}}(\bmu^{ {N}+1})$;
\State Augment RB space $X^\calN_{ {N}+1} = X^\calN_{ {N}} \oplus \{u^\calN(\bmu^{ {N}+1})\}$;
\State $S_{ {N}+1} = S_{ {N}}\bigcup \{\bmu^{ {N}+1}\}$;
\State Set $ {N} \leftarrow  {N}+1$;
\EndWhile
\end{minipage}
}
\vspace{0.05in}
\EndWhile
\end{algorithmic}
\caption{Adaptively constructed STS based classical greedy ${\left(N, X_N^\calN\right) = \mathbb {STS\_CG}}(\Xi_{\rm train}, \varepsilon_{\mathrm{tol}}, C_{M}, K_{\rm damp}, N, X^\calN_N)$}
\label{alg:oerbm}
\end{algorithm}

\section{Hybrid offline enhancements for the greedy algorithms}
\label{sec:multi-greedy}
The unifying theme of the three approaches in the {previous} section is the construction of a small-size subset of the full training set on which the classical greedy algorithm is performed. Algorithm \ref{alg:oerbm} takes full advantage of these information when constructing the small subsets.  Nonetheless, the construction of the surrogate training sets relying on the evaluation of $\Delta_{N}(\bmu)$ for all $\bmu$ in the training set $\Xi_{\rm train}$ is still a bottleneck of Algorithm \ref{alg:oerbm}. 
Algorithms \ref{sen-greedy} and \ref{Jan-greedy} {do not rely on $\Delta_{N}(\bmu)$ and construct the subset by randomly choosing parameters in $\Xi_{\rm train}$}. { However, determining the size of the subset becomes tricky. Too small a subset will degrade the online efficiency, but a very large subset will hinder the offline speed.} 
{Motivated by the thought of combining the strengths of both types of algorithms, we} develop the following two {hybrid} algorithms.
{Indeed, we start with Algorithm \ref{sen-greedy} or \ref{Jan-greedy} while choosing a relatively large initial subset. We then adopt the idea of Algorithm \ref{alg:oerbm} to} build up surrogate training set $\Xi_{\rm sur}$ {by computing the {\it a posteriori} error estimate of the parameters in this subset. 
Next, we details these two hybrid algorithms.
}

\paragraph{Hybrid Training Set Decomposition}
{The idea of Algorithm \ref{sen-greedy} is to decompose the full training set $\Xi_{\rm train}$ into small subsets $\{\Xi_{\rm train, j}\}_{j=1}^{\mathcal{J}}$ and then run standard greedy algorithm on each of them sequentially from $\Xi_{\rm train, 1}$ to $\Xi_{\rm train, \mathcal{J}}$. Since, in practical examples, most snapshots are chosen when performing the classical greedy algorithm on $\Xi_{\rm train, 1}$, $\Xi_{\rm train, 1}$ can not be too small relative to the full training set. Otherwise, it is not rich enough to emulate the full training set. However, large $\Xi_{\rm train, 1}$ increases the offline cost significantly. 
Our first hybrid {method}, {Algorithm \ref{alg:sen-smm} integrating Algorithms \ref{sen-greedy} and \ref{alg:oerbm},} mitigates this dilemma by running Algorithm \ref{alg:oerbm} in each of the subsets generated by the training set decomposition  $\{\Xi_{\rm train, j}\}_{j=1}^{\mathcal{J}}$.}

{That is, g}iven a full training set $\Xi_{\rm train}$, first we construct a training set decomposition according to Algorithm \ref{sen-greedy} to obtain $\{\Xi_{\rm train, j}\}_{j=1}^{\mathcal{J}}$. Then, we perform algorithm \ref{alg:oerbm} on each sample set $\Xi_{\rm train, j}$ sequentially. 
As verified by our numerical results, this simple hybrid approach tends to reduce the greedy algorithm's sensitivity to $n_{\rm tr, small}$, while retaining the ease of picking $C_M$ and $K_{\rm damp}$ as outlined in section \ref{sec:smm-greedy}. 
\begin{algorithm}
\begin{algorithmic}[1]
\State {Determine a partition $\{\Xi_{\rm train, j}\}_{j=1}^\calJ$ for $\Xi_{\rm train}$, e.g. by randomly sampling $n_{\rm train, j}$ points from $\Xi_{\rm train}$ to form $\Xi_{\rm train, j}$ where $|\Xi_{\rm train, j}| = 2^{j-1}n_{\rm tr, small}$ for $j = 1, \dots, \mathcal{J}-1$, and  $\Xi_{\rm train, \calJ}$ completing the partition.}
\State Choose $\bmu^{1} \in \Xi_{\rm train,1}$ at random. Denote $S_{1} = \{ \bmu^{1}\}$, $X^\calN_{1} = {\rm span} \{u^\calN(\bmu^{1}) \}$ and $N = 1$;
\For{$j = 1: \mathcal{J}$}
\State Call $\left( N, X^\calN_N \right) = {\mathbb {STS\_CG}}(\Xi_{\rm train, j}, \varepsilon_{\mathrm{tol}}, C_{M}, K_{\rm damp}, N, X^\calN_N)$
\EndFor
\end{algorithmic}
 \caption{Hybrid training set decomposition-based greedy algorithm $\left(N, X^\calN_N\right)={\mathbb {H\_TSD\_CG}}(\Xi_{\rm train}, \varepsilon_{\mathrm{tol}}, C_{M}, K_{\rm damp})$} 
 \label{alg:sen-smm}
\end{algorithm}

\paragraph{Hybrid Adaptive Enriching}
The randomly-generated sample set $\Xi$ in Algorithm \ref{Jan-greedy}, albeit not covering all of $\Xi_{\rm train}$, keeps updating itself after each round of greedy searching by replacing parameter values whose {\it a posteriori} error estimate falls below $\varepsilon_{\rm tol}$ by those not seen yet. To avoid having too many rounds of replacements, 
$\Xi$ must be rich enough to be representative of $\Xi_{\rm train}$. However, large $\Xi$ adversely impact the algorithm's efficiency. 
The proposed hybrid {strikes} a balance between efficiency and richness. 
Toward that end and as outlined in Algorithm \ref{alg:Jan-smm}, {assimilating Algorithms \ref{Jan-greedy} and \ref{alg:oerbm},} we adaptively construct $\Xi_{\rm sur}$ for the active training set $\Xi$ which is only fully examined when (the much smaller) $\Xi_{\rm sur}$ is sufficiently resolved. This three-level approach enables a faster examination of $\Xi$ than the Adaptive Enriching algorithm \ref{Jan-greedy} which, in turn, allows $\Xi$ to be larger thus more representative of $\Xi_{\rm train}$.
\begin{algorithm}
\begin{algorithmic}[1]
\State $N_{safe} = {\rm ceil}(|\Xi_{\rm train}|/M_{\rm sample})$;
\State Randomly generate an initial training set $\Xi$ with $M_{\rm sample}$ parameters from $\Xi_{\rm train}$;
\State Choose an initial parameter $\bmu^{1} \in \Xi_{\rm train}$.  Set $S_{1} = \{ \bmu^{1}\}$, $X^\calN_{1} = {\rm span} \{u^\calN(\bmu^{1}) \}$ and $N = 1$;
\State Set ${\rm safe} =0, \ell = 1$, $\varepsilon = 2\varepsilon_{\mathrm{tol}}$ and $r = n_{\rm train}$;
\While{(${E_{\ell}} > \varepsilon_{\mathrm{tol}} $ or ${\rm safe} \le N_{\rm safe}$) and ($r > 0$)}
{\State
\begin{minipage}{0.03\textwidth}
\rotatebox{90}{One-step greedy}
\end{minipage} 
\begin{minipage}{0.03\textwidth}
\rotatebox{90}{scan on $\Xi$}
\end{minipage} 
\hspace{4mm}
%one-step greedy starts
\framebox[0.8 \textwidth]
{\begin{minipage}{0.8 \textwidth}
\For{each $\boldsymbol{\mu} \in \Xi$}
\State Obtain RBM solution $u^{\mathcal{N}}_{ {N}}(\boldsymbol{\mu}) \in X^\calN_{ {N}}$ and error estimate ${\Delta_{ {N}}}(\boldsymbol{\mu})$;
\EndFor
\vspace{0.05in}
\State $\bmu^{ {N}+1} = \underset{\bmu \in \Xi}{\argmax} \Delta_{ {N}}(\bmu)$, $\varepsilon = \Delta_{ {N}}(\bmu^{ {N}+1})$, $E_{\ell} = \varepsilon$;
\State Augment RB space $X^\calN_{ {N}+1} = X^\calN_{ {N}} \oplus \{u^\calN(\bmu^{ {N}+1})\}$;
\State $S_{ {N}+1} = S_{ {N}}\bigcup \{\bmu^{ {N}+1}\}$;
\State Set $ {N} \leftarrow  {N}+1$;
\end{minipage}
}
%one-step greedy ends
\vspace{0.05in}
\State \hspace{-0.0\textwidth} Construct STS $\Xi_{\rm Sur}$ based on $\{(u_{ {N}-1}^\calN(\bmu), \Delta_{ {N}-1}(\bmu)): \bmu \in \Xi\}$ with SMM;
\vspace{0.05in}
\State 
\begin{minipage}{0.03\textwidth}
\rotatebox{90}{Multi-step greedy}
\end{minipage} 
\begin{minipage}{0.03\textwidth}
\rotatebox{90}{scan on $\Xi_{\rm Sur}$}
\end{minipage} 
\hspace{4mm}
%multi-step greedy starts
\framebox[0.8 \textwidth]
{\begin{minipage}{0.8 \textwidth}
\While {$(\varepsilon > \varepsilon_{\mathrm{tol}})$ and {$(\frac{\varepsilon}{E_{\ell}} >  \, \frac{1}{K_{\rm damp} \times (\ell + 1)})$}}
\vspace{0.05in}
\For{each $\boldsymbol{\mu} \in \Xi_{\rm Sur}$}
\State Obtain RBM solution $u^{\mathcal{N}}_{ {N}}(\boldsymbol{\mu}) \in X^\calN_{ {N}}$ and error estimate ${\Delta_{ {N}}}(\boldsymbol{\mu})$;
\EndFor
\vspace{0.05in}
\State $\bmu^{ {N}+1} = \underset{\bmu \in \Xi_{\rm Sur}}{\argmax} \Delta_{ {N}}(\bmu)$, $\varepsilon = \Delta_{ {N}}(\bmu^{ {N}+1})$;
\State Augment RB space $X^\calN_{ {N}+1} = X^\calN_{ {N}} \oplus \{u^\calN(\bmu^{ {N}+1})\}$;
\State $S_{ {N}+1} = S_{ {N}}\bigcup \{\bmu^{ {N}+1}\}$;
\State Set $ {N} \leftarrow  {N}+1$;
\EndWhile
\end{minipage}
}
}
\State Truncate $\Xi$ by $\Xi_{<\varepsilon} = \{ \bmu \in \Xi: \Delta_N(\bmu) < \varepsilon_{\mathrm{tol}}\}$;
\State Truncate $\Xi_{\rm train}$ by $\Xi_{<\varepsilon}$, and set $r = r -  |\Xi_{<\varepsilon}|$;
\If{$|\Xi_{<\varepsilon}| = M_{\rm sample}$}
\State Set ${\rm safe} =  {\rm safe} +1$, and $\ell={1}$;
{
\Else
\State $\ell = \ell + 1;$}
\EndIf
\State Randomly choose $M_{\rm sample} - |\Xi|$ parameters from $\Xi_{\rm train}$ for addition to $\Xi$;
\EndWhile
\end{algorithmic}
 \caption{Hybrid adaptive enriching greedy algorithm ${\mathbb {H\_AE\_CG}}(\Xi_{\rm train}, \varepsilon_{\mathrm{tol}}, C_{M}, K_{\rm damp},  {M_{\rm sample}})$}
  \label{alg:Jan-smm}
\end{algorithm}

{
\begin{remark}
As pointed out in Sections \ref{sec:senRB} and \ref{sec: AE} and shown in the numerical results section, the performance of the non-hybrid Algorithms  \ref{sen-greedy} and \ref{Jan-greedy} exhibit high sensitivity with respect to the tuning parameters. 
A direct consequence is that selecting an appropriate tuning parameter is not easy. As will be shown in the numerical results section, their hybrid version Algorithms \ref{alg:sen-smm} and \ref{alg:Jan-smm} 
not only come with computational acceleration but also are capable of significantly reducing the sensitivity of the tuning parameters thus making the methods more usable.
\end{remark}
}

\section{Numerical Tests}
\label{sec: numerics}

In this section, we test and compare these six greedy algorithms, namely classical greedy, three existing improvements reviewed in section \ref{sec:background}, and our newly designed two hybrids in section \ref{sec:multi-greedy}. For the sake of clarity, we list these methods in one place, Table \ref{tab:notation-numerics}. 
\begin{table}
  \begin{center}
  \resizebox{0.8\textwidth}{!}{
    \renewcommand{\tabcolsep}{0.6cm}
    \renewcommand{\arraystretch}{1.3}
 {
\begin{tabular}{@{}lp{0.65\textwidth}@{}}
      \toprule    
Abbrv & Full name \\
\midrule
${\mathbb {CG}}$ & Classical greedy algorithm\\   
${\mathbb {TSD\_CG}}$ & Classical greedy enhanced by Training Set Decomposition\\
${\mathbb {AE\_CG}}$ & Classical greedy enhanced by Adaptive Enriching\\
${\mathbb {STS\_CG}}$ & Classical greedy enhanced by adaptively constructed Surrogate Training Set\\
${\mathbb {H\_TSD\_CG}}$ & Hybrid training set decomposition based greedy algorithm\\
${\mathbb {H\_AE\_CG}}$ & Hybrid adaptive enriching greedy algorithm\\ 
 \bottomrule
    \end{tabular}
}
\renewcommand{\arraystretch}{1}
    \renewcommand{\tabcolsep}{12pt}
}
  \end{center}
  \caption{Abbreviation of the greedy algorithms used in this section.}\label{tab:notation-numerics}
\end{table}

Two types of examples will be presented to demonstrate the efficiency enhancement of our proposed approaches without sacrificing the quality of the reduced bases. The corresponding results are presented in the subsections below. The CPU times reported in this paper refer to computations performed on a workstation with 3.1 GHz Intel Core i5 processor and 16GB memory, in the MATLAB environment adopting redbKIT library \cite{redbKIT, Quarteroni2015} {and RBmatlab package \cite{RBmatlab, RBmatlabPaper}\footnote{Available for download at {\texttt http://www.ians.uni-stuttgart.de/MoRePaS/software/index.html}}}.

\subsection{{Thermal Block problem with nine parameters}}
\label{sec:diff}
{We start with the classical thermal block problem \cite{Rozza2008, Quarteroni2015, sen2008reduced, jiang2017offline}. 
\begin{equation}\label{eq:diff_prob_2}
\begin{cases}
- \nabla . (a(x, \bmu) \nabla u(x, \bmu) ) = f \quad \text{on}  \quad \Omega,\\
u(x,\bmu) = g_{D} \quad \text{on} \quad \Gamma_{D}, \\
\frac{\partial u}{\partial n} = g_{N} \quad \text{on} \quad \Gamma_{N}.
\end{cases}
\end{equation}
Here, $\Omega = [0, 1] \times [0, 1]$ is partitioned into $9$ blocks $\bigcup_{i=1}^9 B_i = \Omega$, $\Gamma_{D}$ is the top boundary, and $\Gamma_{N}  = \partial \Omega   \setminus \Gamma_{D}$. The heat conductivities on these blocks constitute the parameters of our test:
\begin{center}
\begin{tabular}{|c|c|c|}
\multicolumn{3}{c}{$\Gamma_{D}$} \\
\hline
$\mu_{7} (B_{7})$ & $\mu_{8} (B_{8})$ & $\mu_{9} (B_{9})$ \\
\hline
$\mu_{4} (B_{4})$ & $\mu_{5} (B_{5})$ & $\mu_{6} (B_{6})$ \\
\hline
$\mu_{1} (B_{1})$ & $\mu_{2} (B_{2})$& $\mu_{3} (B_{3})$ \\
\hline
\multicolumn{3}{c}{$\Gamma_{base}$} 
\end{tabular}
\end{center}
That is, the diffusion coefficient $a(x, \bmu) = \mu_{i} ~\text{when} ~x \in B_{i}$. The parameter vector is thus given by $\bmu = (\mu_{1}, \mu_{2}, \dots, \mu_{9})$ whose domain is chosen as $\mathcal{D} = [0.1, 10]^9$ for our test. We take as the right hand side $f = 0$. The Dirichlet boundary is set to be homogeneous, i.e. $g_{D} = 0$. The Neumann data is set to simulate heat influx only at the bottom, i.e. $g_{N} = 1$ on the bottom boundary $\Gamma_{\rm base}$ and $g_{N} = 0$ at the two sides. The output of interest is defined to be the average temperature over $\Gamma_{\rm base}$
\begin{equation}\label{eq: output}
s(\bmu) = \int_{\Gamma_{\rm base}} u(x,\bmu)dx.
\end{equation}
The truth approximation is obtained by a finite element solver \cite{RBmatlab, RBmatlabPaper} with total degrees of freedom $\mathcal{N} = 361$. A sufficient number of samples ($N_{{\rm train}} = {4\times 10^6}$) are randomly drawn from the parameter domain $\mathcal{D}$. The testing set contains another $1,000$ random samples in $\mathcal{D}$. The error is measured in $H^1$ norm the largest of which is recorded to show the worst case scenario. Moreover, to gain an understanding of the probabilistic nature of the algorithms, we run each of the ${\mathbb {TSD\_CG}}, {\mathbb {AE\_CG}},{\mathbb {H\_TSD\_CG}}, {\mathbb {H\_AE\_CG}}$ $5$ times for every fixed tuning parameter and tabulate the corresponding ranges.}

\begin{figure}[htbp]
\centering
\includegraphics[width = 0.7\textwidth]{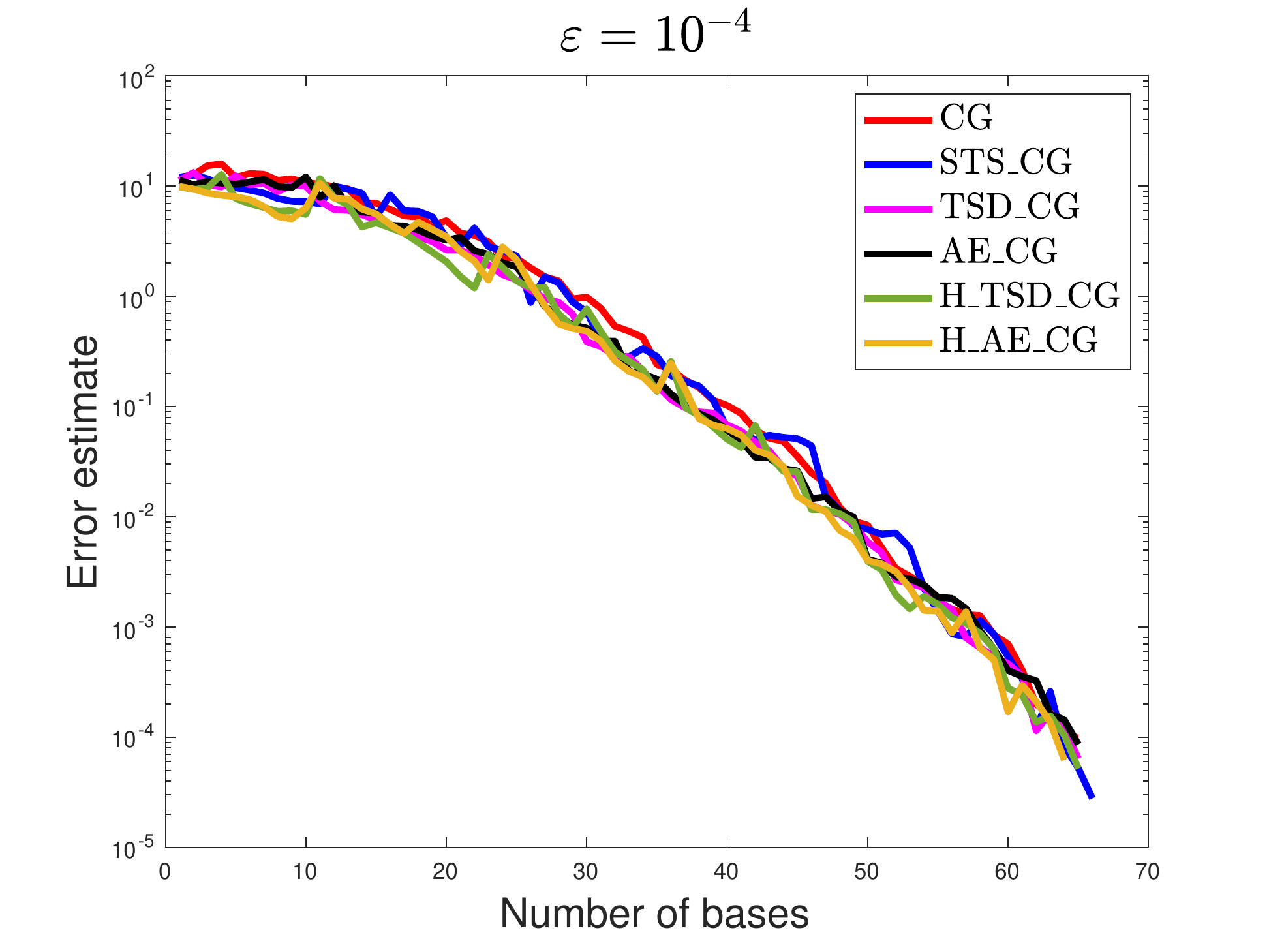}
\caption{{Convergence of a posteriori error estimate for the thermal block problem.}}
\label{fig:ex1_err}
\end{figure}

First, we fix $K_{\rm damp} =  {20}$, $C_M = {20}$, $N_{\rm sample} = {2\times 10^5}$ and $n_{\rm tr,small} = {2\times 10^5}$ and test {3} different tolerances. Table \ref{eff_num} shows {the accuracy and} numbers of bases at convergence for each method, and its speedup factor in comparison to the classical greedy algorithm. For relative time, defined as the corresponding running time scaled by the running time of classical greedy algorithm, we observe that the worst case scenario of {the two newly proposed approaches, ${\mathbb {H\_TSD\_CG}}$ {and ${\mathbb {H\_AE\_CG}}$},} is better than the best case scenario of any other {approaches}. 
Moreover, {the efficiency is increasing as tolerance gets smaller. Therefore, the alleviation is more pronounced for high-dimensional problem or when high accuracy is desired. In addition, with the same tolerance, the required number of bases are very similar for all these methods. It indicates that offline speedup of both ${\mathbb {H\_TSD\_CG}}$ and ${\mathbb {H\_AE\_CG}}$ are not coming at the cost of online efficiency. In terms of accuracy, the worst case (i.e. upper bounds of errors) for ${\mathbb {H\_TSD\_CG}}$ and ${\mathbb {H\_AE\_CG}}$ are in the same order of magnitude as the $\mathbb {CG}$. Figure \ref{fig:ex1_err} illustrates that the a posteriori error estimate is converging exponentially for both ${\mathbb {H\_TSD\_CG}}$ and ${\mathbb {H\_AE\_CG}}$, in the same fashion as other approaches including the canonical greedy algorithm. This shows that our proposed hybrid algorithms do not appear to suffer accuracy degradation for this example.}

\begin{table}[h!]
  \centering

  \medskip

  \begin{tabularx}{\linewidth}{ X X X X X X X}
    \multicolumn{6}{c}{(a) Efficiency in terms of relative computation time} \\
    \toprule
    $\varepsilon_{\rm tol}$&${\mathbb {CG}}$&${\mathbb {STS\_CG}}$&${\mathbb {TSD\_CG}}$&${\mathbb {AE\_CG}}$&${\mathbb {H\_TSD\_CG}}$&${\mathbb {H\_AE\_CG}}$\\
    \midrule
   $10^{-2} $& 1 & {0.103}  & ${[0.055, 0.057]}$ & ${[0.061, 0.068]}$ & ${[0.035,    0.036]}$ &${[0.033,    0.034]}$\\
    \midrule
    $10^{-3}$ &  1 & {0.106}& ${[0.056, 0.063]}$ & ${[0.060, 0.064]}$ &${[0.035, 0.036]}$& ${[0.031, 0.032]}$\\
    \midrule
    $ 10^{-4}$ &  1 & {0.097} & ${[0.052, 0.053]}$  & ${[0.058, 0.063]}$ &${[0.028, 0.029]}$& ${[0.028, 0.029]}$\\
    \bottomrule
  \end{tabularx}

  \bigskip
  \begin{tabularx}{\linewidth}{ X X X X X X X}
    \multicolumn{6}{c}{(b) Number of bases at convergence} \\
    \toprule
    $\varepsilon_{\rm tol}$&${\mathbb {CG}}$&${\mathbb {STS\_CG}}$&${\mathbb {TSD\_CG}}$&${\mathbb {AE\_CG}}$&${\mathbb {H\_TSD\_CG}}$& ${\mathbb {H\_AE\_CG}}$\\
    \midrule
     $10^{-2} $& {50} &{52} & ${[50, 51]}$  & ${[50, 51]}$  & ${[51, 51]}$&${[50, 52]}$\\
    \midrule
    $10^{-3}$ & {60} & {60} &  ${[58, 60]}$ &${[59, 60]}$ &${[59, 59]}$&${[59, 61]}$\\
    \midrule
    $10^{-4}$ & {65} & {66} &  ${[64, 66]}$  &${[65, 65]}$ &${[64, 66]}$&${[65, 66]}$\\
    \bottomrule
  \end{tabularx}

  \bigskip
  
  \begin{tabularx}{\linewidth}{ X X X X X X X}
    \multicolumn{6}{c}{ {(c) Accuracy($H^1$ norm): (in $10^{-4}$)}} \\
    \toprule
    $\varepsilon_{\rm tol}$&${\mathbb {CG}}$&${\mathbb {STS\_CG}}$&${\mathbb {TSD\_CG}}$&${\mathbb {AE\_CG}}$&${\mathbb {H\_TSD\_CG}}$& ${\mathbb {H\_AE\_CG}}$\\
    \midrule
     $10^{-2} $& {$1.34$} &{$0.84$} & ${[1.11, 2.93]}$  & ${[0.74, 3.11]}$  & ${[1.04, 1.40]}$&${[1.46, 2.94]}$\\
    \midrule
    $10^{-3}$ & {0.21} & {0.09} &  ${[0.08, 0.44]}$ &${[0.11, 0.31]}$ &${[0.13, 0.28]}$&${[0.05, 0.22]}$\\
    \midrule
    $10^{-4}$ & {0.009} & {0.008} &  ${[0.008, 0.01]}$  &${[0.006, 0.01]}$ &${[0.008, 0.01]}$&${[0.006, 0.008]}$\\
    \bottomrule
  \end{tabularx} 
  
    \caption{{Results for the nine-dimensional thermal block problem: (a) Comparison of {speedup}, (b) size of the resulting reduced basis space, and (c) accuracy.}}
  \label{eff_num}
\end{table}

Next, we consider the impact of tweaking the tuning parameters.  Toward that end, for ${\varepsilon_{\rm tol} = 10^{-3}}$, we choose $(K_{\rm damp}, C_M) = ( 5, {5}), ( 10, {10}), ( 20, {20})$, $N_{\rm sample} = {1\times 10^5, 2\times 10^5, 4\times 10^5}$ and $n_{\rm tr,small} =  {1\times 10^5, 2\times 10^5, 4\times 10^5}$ and record the relative time, the numbers of bases at convergence { and the RBM error}. The results for  {comparing} ${\mathbb {STS\_CG}}$, ${\mathbb {TSD\_CG}}$, ${\mathbb {AE\_CG}}$, ${\mathbb {H\_TSD\_CG}}$ {and ${\mathbb {H\_AE\_CG}}$} are shown in {Table \ref{tab:tun-eff}, Table \ref{tab:tun-basis} and Table \ref{tab:tun-error}.} 
For fixed tolerance, it is evident that the performance of the newly proposed ${\mathbb {H\_TSD\_CG}}$ {and ${\mathbb {H\_AE\_CG}}$} is very stable for various selection of tuning parameters in speed, size of RB space {and accuracy. For relative computational time, the largest of the relative time variations (defined as $(t_{max}-t_{min})/t_{min}$) of these two methods is $31\%$.}  On the other hand, {the efficiency of } ${\mathbb {TSD\_CG}}$ and ${\mathbb {AE\_CG}}$ {is severely sensitive of the} tuning {parameters}. {In fact, we notice that the setup with the worst tuning parameter settings is half as slow as the that of the optimal setting. Note that although $\mathbb {STS\_CG}$ is as stable as ${\mathbb {H\_TSD\_CG}}$ and ${\mathbb {H\_AE\_CG}}$,  its relative computational time is more than three times higher. Although ${\mathbb {H\_TSD\_CG}}$ and ${\mathbb {H\_AE\_CG}}$ are not deterministic, the size of RB space and the accuracy of the resulting reduced solver are both stable with the accuracy being comparable to  that of the classical greedy algorithm. 
}

{Based on these two experiments, we clearly see that our newly proposed ${\mathbb {H\_TSD\_CG}}$ and ${\mathbb {H\_AE\_CG}}$ improve} both efficiency and stability. {Furthermore, the} reduction in tuning parameter sensitivity makes the hybrids particularly attractive.

\begin{table}
\begin{tabularx}{\linewidth}{ X X X X }
    \multicolumn{4}{c}{(a) ${\mathbb {STS\_CG}}$} \\
    \toprule
    $(K_{\rm damp},C_M)$ & ( {5}, {5}) & ( {10}, {10}) & ( {20}, {20})\\
    \midrule
     & {0.129} & {0.127} & {0.106} \\
    \bottomrule
  \end{tabularx}
  
  \bigskip
  
  \begin{tabularx}{\linewidth}{ X X X X }
    \multicolumn{4}{c}{(b) ${\mathbb {TSD\_CG}}$} \\
    \toprule
    $ n_{\rm tr,small}$ &  ${1 \times 10^5}$ & ${2\times 10^5}$ & ${4\times 10^5}$\\
    \midrule
     & ${[0.041, 0.047]}$ & ${[0.056, 0.063]}$ &  ${[0.089, 0.10]}$  \\
    \bottomrule
  \end{tabularx}
  
  \bigskip
  
  \begin{tabularx}{\linewidth}{ X X X X }
    \multicolumn{4}{c}{(c) ${\mathbb {AE\_CG}}$} \\
    \toprule
    $M_{\rm sample}$ & ${1 \times 10^5}$ & ${2\times 10^5}$ & ${4\times 10^5}$\\
    \midrule
     & ${[0.044,    0.045]}$ & ${[0.060,    0.064]}$ &  ${[0.086,    0.097]}$\\
    \bottomrule
  \end{tabularx}
  
  \bigskip
  
  \begin{tabularx}{\linewidth}{ X X X X X X}
    \multicolumn{6}{c}{(d) ${\mathbb {H\_TSD\_CG}}$} \\
    \toprule
    \multicolumn{3}{c}{$(K_{\rm damp}, C_M) / n_{\rm tr,small}$} & ${1 \times 10^5}$ & ${2\times 10^5}$ & ${4\times 10^5}$\\
    \midrule 
    \multicolumn{3}{c}{ $( {5}, {5})$} & ${[0.030,0.031]}$ & ${[0.036,    0.037]}$ & ${[0.037,    0.038]}$  \\
      \multicolumn{3}{c}{$( {10}, {10})$} &${[0.030,0.030]}$ &${[0.036,    0.037]}$ &${[0.036,    0.037]}$ \\
      \multicolumn{3}{c}{$( {20}, {20})$} &${[0.029,0.030]}$ &${[0.035,    0.036]}$ &${[0.035,    0.036]}$  \\
    \bottomrule
  \end{tabularx}  
  
  \bigskip
  
  \begin{tabularx}{\linewidth}{ X X X X X X}
    \multicolumn{6}{c}{(e) ${\mathbb {H\_AE\_CG}}$} \\
    \toprule
    \multicolumn{3}{c}{$(K_{\rm damp}, C_M) / n_{\rm tr,small}$} & ${1 \times 10^5}$ & ${2\times 10^5}$ & ${4\times 10^5}$\\
    \midrule 
    \multicolumn{3}{c}{ $( {5}, {5})$} & ${[0.030,    0.031]}$ & ${[0.032,    0.034]}$ & ${[ 0.033,    0.038]}$  \\
      \multicolumn{3}{c}{$( {10}, {10})$} &${[0.029,    0.030]}$ &${[0.031,    0.032]}$ &${[0.035,    0.036]}$ \\
      \multicolumn{3}{c}{$( {20}, {20})$} &${[0.029,    0.030]}$ &${[0.031,    0.032]}$ &${[0.034,    0.035]}$  \\
    \bottomrule
  \end{tabularx} 
 \caption{Tuning parameter study for the {thermal block problem}. Comparison of {relative efficiency for various enhanced greedy algorithms.}}
  \label{tab:tun-eff}
\end{table}

\begin{table}
\begin{tabularx}{\linewidth}{ X X X X }
    \multicolumn{4}{c}{(a) ${\mathbb {STS\_CG}}$} \\
    \toprule
    $(K_{\rm damp},C_M)$ & ( {5}, {5}) & ( {10}, {10}) & ( {20}, {20})\\
    \midrule
     & {60} & {60} & {60} \\
    \bottomrule
  \end{tabularx}
  
  \bigskip
  
  \begin{tabularx}{\linewidth}{ X X X X }
    \multicolumn{4}{c}{(b) ${\mathbb {TSD\_CG}}$} \\
    \toprule
    $n_{\rm tr,small}$ & ${1 \times 10^5}$ & ${2\times 10^5}$ & ${4\times 10^5}$\\
    \midrule
     & ${[59,   60]}$ & ${[58,   60]}$ & ${[59, 59]}$ \\
    \bottomrule
  \end{tabularx}
  
  \bigskip
  
  \begin{tabularx}{\linewidth}{ X X X X }
    \multicolumn{4}{c}{(c) ${\mathbb {AE\_CG}}$} \\
    \toprule
    $M_{\rm sample}$ & ${1 \times 10^5}$ & ${2\times 10^5}$ & ${4\times 10^5}$\\
    \midrule
     & ${[59,   60]}$ & ${[59,   60]}$ &  ${[57, 58]}$\\
    \bottomrule
  \end{tabularx}
  
  \bigskip
  
  \begin{tabularx}{\linewidth}{ X X X X X X X}
    \multicolumn{6}{c}{(d) ${\mathbb {H\_TSD\_CG}}$} \\
    \toprule
    \multicolumn{3}{c}{$(K_{\rm damp}, C_M) / n_{\rm tr,small}$} & ${1 \times 10^5}$ & ${2\times 10^5}$ & ${4\times 10^5}$\\
    \midrule 
    \multicolumn{3}{c}{ $( {5}, {5})$} &${[59,   60]}$ &${[58,   60]}$ &${[59,   61]}$  \\
      \multicolumn{3}{c}{$( {10}, {10})$} &${[59,   61]}$ &${[59,   62]}$ &${[59,   60]}$ \\
      \multicolumn{3}{c}{$( {20}, {20})$} &${[ 59,   60]}$ &${[59, 59]}$ &${[59,   60]}$  \\
       \bottomrule
  \end{tabularx}  
  
  \bigskip
  
  \begin{tabularx}{\linewidth}{ X X X X X X X}
    \multicolumn{6}{c}{(e) ${\mathbb {H\_AE\_CG}}$} \\
    \toprule
    \multicolumn{3}{c}{$(K_{\rm damp}, C_M) / n_{\rm tr,small}$} & ${1 \times 10^5}$ & ${2\times 10^5}$ & ${4\times 10^5}$\\
    \midrule 
    \multicolumn{3}{c}{ $( {5}, {5})$} &${[59,   62]}$ &${[59,   61]}$ &${[59,   60]}$  \\
      \multicolumn{3}{c}{$( {10}, {10})$} &${[59,   61]}$ &${[59,   61]}$ &${[59,   60]}$ \\
      \multicolumn{3}{c}{$( {20}, {20})$} &${[59,   61]}$ &${[59,   61]}$ &${[59,   60]}$  \\
       \bottomrule
  \end{tabularx}  
\caption{Tuning parameter study for the {thermal block problem}. Comparison of number of bases at convergence {for various enhanced greedy algorithms.}}
  \label{tab:tun-basis}
\end{table}

\begin{table}
\begin{tabularx}{\linewidth}{ X X X X }
    \multicolumn{4}{c}{(a) ${\mathbb {STS\_CG}}$} \\
    \toprule
    $(K_{\rm damp},C_M)$ & ( {5}, {5}) & ( {10}, {10}) & ( {20}, {20})\\
    \midrule
     & {0.12} & {0.05} & {0.09} \\
    \bottomrule
  \end{tabularx}
  
  \bigskip
  
  \begin{tabularx}{\linewidth}{ X X X X }
    \multicolumn{4}{c}{(b) ${\mathbb {TSD\_CG}}$} \\
    \toprule
    $n_{\rm tr,small}$ & ${1 \times 10^5}$ & ${2\times 10^5}$ & ${4\times 10^5}$\\
    \midrule
     & ${[0.07,   0.17]}$ & ${[0.08,   0.4]}$ & ${[0.11, 0.27]}$ \\
    \bottomrule
  \end{tabularx}
  
  \bigskip
  
  \begin{tabularx}{\linewidth}{ X X X X }
    \multicolumn{4}{c}{(c) ${\mathbb {AE\_CG}}$} \\
    \toprule
    $M_{\rm sample}$ & ${1 \times 10^5}$ & ${2\times 10^5}$ & ${4\times 10^5}$\\
    \midrule
     & ${[0.07,   0.24]}$ & ${[0.11,   0.31]}$ &  ${[0.25, 0.51]}$\\
    \bottomrule
  \end{tabularx}
  
  \bigskip
  
  \begin{tabularx}{\linewidth}{ X X X X X X X}
    \multicolumn{6}{c}{(d) ${\mathbb {H\_TSD\_CG}}$} \\
    \toprule
    \multicolumn{3}{c}{$(K_{\rm damp}, C_M) / n_{\rm tr,small}$} & ${1 \times 10^5}$ & ${2\times 10^5}$ & ${4\times 10^5}$\\
    \midrule 
    \multicolumn{3}{c}{ $( {5}, {5})$} &${[0.08,   0.21]}$ &${[0.11,   0.24]}$ &${[0.07,   0.37]}$  \\
      \multicolumn{3}{c}{$( {10}, {10})$} &${[0.05,   0.18]}$ &${[0.04,   0.18]}$ &${[0.10,   0.21]}$ \\
      \multicolumn{3}{c}{$( {20}, {20})$} &${[ 0.11,   0.28]}$ &${[0.13, 0.28]}$ &${[0.19,   0.51]}$  \\
       \bottomrule
  \end{tabularx}  
  
  \bigskip
  
  \begin{tabularx}{\linewidth}{ X X X X X X X}
    \multicolumn{6}{c}{(e) ${\mathbb {H\_AE\_CG}}$} \\
    \toprule
    \multicolumn{3}{c}{$(K_{\rm damp}, C_M) / n_{\rm tr,small}$} & ${1 \times 10^5}$ & ${2\times 10^5}$ & ${4\times 10^5}$\\
    \midrule 
    \multicolumn{3}{c}{ $( {5}, {5})$} &${[0.03,   0.22]}$ &${[0.10,   0.24]}$ &${[0.11,   0.22]}$  \\
      \multicolumn{3}{c}{$( {10}, {10})$} &${[0.09,   0.18]}$ &${[0.08,   0.25]}$ &${[0.07,   0.28]}$ \\
      \multicolumn{3}{c}{$( {20}, {20})$} &${[0.07,   0.19]}$ &${[0.05,   0.22]}$ &${[0.08,   0.26]}$  \\
       \bottomrule
  \end{tabularx}  
\caption{{Tuning parameter study for the thermal block problem. Comparison of algorithmic accuracy for various enhanced greedy algorithms. $H^1$ norm (in $10^{-4}$).}}
  \label{tab:tun-error}
\end{table}

\subsection{Helmholtz equation on a parametrized domain}
\label{sec:helm}

\begin{figure}[htbp]
\centering
\includegraphics[width=0.7\textwidth]{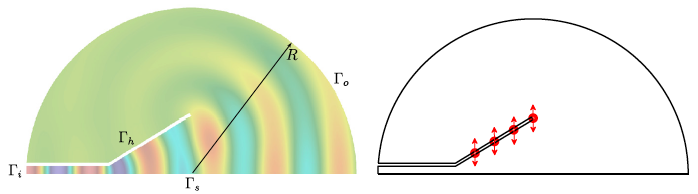}
\caption{Left: the acoustic horn domain and boundaries (background coloring given by $\Re(p)$ for $f = 900 $Hz). Right: RBF control points (red circles) whose vertical displacement is treated as a parameter.  [Figure taken from \cite{negri2015efficient}].}
\label{fig:hornsetup}
\end{figure}

{For our second test, we turn to a more challenging non-coercive and nonaffine problem which is parameterized by the geometric configuration of the system. Indeed, w}e consider the propagation of {a} pressure wave $P(\bf{x},t)$ into the acoustic horn  illustrated in Figure \ref{fig:hornsetup}, the same example considered in \cite{negri2015efficient}. {Assuming that the} waves {are} time harmonic, the acoustic pressure $P$ can be separated as $P(x, t) = \Re(p(x)\exp^{iwt})$ where the complex amplitude $p(x)$ satisfies the following Helmholtz equation \cite{kasolis2012fixed, negri2015efficient}:

\begin{align}
\Delta p + \kappa^2p                                            &= 0                &   \text{in}\quad  \Omega \nonumber\\
(i\kappa + \frac{1}{2R})p + \nabla p \cdot {\bf{n}} & = 0               &    \text{on} \quad \Gamma_{o}\nonumber\\
i\kappa p + \nabla p \cdot {\bf{n}}                         &= 2i\kappa A  &    \text{on} \quad \Gamma_{i} \nonumber\\
\nabla p \cdot {\bf{n}}                                            &= 0               &   \text{on} \quad \Gamma_{h}\cup \Gamma_{s} = \Gamma_{n},
\label{eq:hel_prob}
\end{align}
where $\kappa = w/c$ is the wave number, $w = 2\pi f$ the angular frequency and $c = 340 \frac{\rm cm}{\rm s}$ the speed of sound. A radiation condition is prescribed on the boundary $\Gamma_{i}$ imposing an inner-going wave with amplitude $A=1$ and absorbing the outer-going planar waves. A Neumann boundary condition is applied on the walls $\Gamma_{h}$ of the device as well as on the symmetry boundary $\Gamma_{s}$. Finally, an absorbing condition is placed on the far-field boundary $\Gamma_{o}$ with radius $R = 1$.

We consider up to five parameters. The first is the frequency $f$.  The other four describe the shape of the horn, representing the vertical displacement of the RBF control points \cite{buhmann2003radial} in Figure \ref{fig:hornsetup}. As a result,  we have the five-dimensional parameter vector  $\bmu = [f \quad \bmu_{g}]$. The output of interest is the index of reflection intensity (IRI)\cite{bangtsson2003shape} defined as 
$$J(\mu) = \left| \frac{1}{\Gamma_{i}}\int_{\Gamma_{i}}p(\mu)d\Gamma -1 \right|, $$
which measures the transmission efficiency of the device. {Since the pressure $p$ depends on the design of the horn and the frequency, we parametrized the reference domain $\Omega$ and define it as $\Omega(\bmu_g)$, where the details of the parametrized construction can be found in \cite{manzoni2012model,quarteroni2015reduced}. Given $\bmu \in \mathcal{D}$, the weak form of problem (\ref{eq:hel_prob}) reads: find $p \in H^1(\Omega(\bmu_g))$ such that
\begin{equation}
a(p, v, \bmu) = g(v,\bmu) \quad \forall v \in H^1(\Omega(\bmu_g)),
\label{weak-helm}
\end{equation}
where 
\begin{align}
\label{eq:apv}
a(p, v, \bmu) &= \int_{\Omega(\bmu_g)} (\nabla p \cdot \nabla \bar{v} - \kappa^2p\bar{v}) d\Omega + i\kappa \int_{\Gamma_o \bigcup \Gamma_i} p \bar{v} d\Gamma + \frac{1}{2R}\int_{\Gamma_o}p \bar{v} d\Gamma, \\
\label{eq:gv}
g(v, \bmu) &= 2i\kappa A \int_{\Gamma_i}\bar{v} d\Gamma.
\end{align}
We adopt the same implementation as \cite{negri2015efficient}, i.e. a conforming triangular finite element method discretizing problem (\ref{weak-helm}) by approximating $H^1(\Omega(\bmu_g))$ with a set of $\mathcal{N}$ piecewise polynomial nodal basis functions $\{\phi_i\}_{i=1}^{\mathcal{N}}$. Then the {FOM} solver ends up with solving the following large linear system:
\begin{equation}
\label{eq:largelinear}
\mathcal{A}(\bmu)\boldsymbol{p} = \boldsymbol{g}(\bmu), 
\end{equation}
where 
\begin{align}
\mathcal{A}_{ij}(\bmu) &= a(\phi_j, \phi_i, \bmu), \\
\boldsymbol{g}_i(\bmu) &= g(\phi_i, \bmu), \quad 1 \le i,j \le \mathcal{N}.
\end{align}
}

{Many scenarios of this problem, e.g. optimizing the shape for the horn's transmission efficiency, require solving the large system (\ref{eq:largelinear}) for many different configurations. 
However, even one such {solution} can be computationally intensive {to obtain} due to the nonaffinity of the problem introduced by the geometric parameterization. Indeed, empirical interpolation method \cite{Barrault_Nguyen_Maday_Patera, Grepl_Maday_Nguyen_Patera, ChaturantabutSorensen2010} is usually used to achieve the online independence of the full model degrees of freedom by approximating the non-affine parametrized PDE with a linear combination of $Q_a$ affine terms. Let us consider for instance the matrix corresponding to the first term of \eqref{eq:apv},  $\int_{\Omega(\bmu_g)} \nabla p \cdot \nabla \bar{v} d \Omega$.
\[
\mathcal{K}_{ij}(\bmu_g) = \int_{\Omega(\bmu_g)} \nabla \phi_j \cdot \nabla \bar{\phi_i} d \Omega = \int_{\widetilde{\Omega}} \mathcal{T}(\bx; \bmu_g) \nabla \phi_j \cdot \nabla \bar{\phi_i} d \Omega.
\]
Here $\widetilde{\Omega}$ is (fixed) reference configuration that is mapped to $\Omega(\bmu_g)$ through a parametric map
\[
{\mathbf F}: \widetilde{\Omega} \times \calD \rightarrow {\mathbb R}^2 \quad \mbox{ such that }\quad {\mathbf F}(\widetilde{\Omega}; \bmu_g) = \Omega(\bmu_g).
\]
 The $2 \times 2$ matrix $\mathcal{T}(\bx; \bmu_g)$ can then be evaluated 
\[
\mathcal{T}(\bx; \bmu_g) = \left(\nabla_x {\mathbf F}(x, \bmu_g)\right)^{-1} \left(\nabla_x {\mathbf F}(x, \bmu_g)\right)^{-T} \left| \nabla_x {\mathbf F}(x, \bmu_g)\right|
\]
}

{The EIM {approximates} each component of $\mathcal{T}(x; \bmu_g)$, for any $\bmu_g$, by its projection onto a low-dimensional space spanned by some well-chosen instances of $\mathcal{T}(x; \bmu_g)$, $\{\mathcal{T}(x; \bmu_g): \bmu_g = \bmu_g^1, \dots, \bmu_g^{Q_{k \ell}}\}$:
\[
(\mathcal{T}(x; \bmu_g))_{k,\ell} = \sum_{q = 1}^{Q_{k \ell}} (\mathcal{T}(x; \bmu_g^q))_{k,\ell} \theta(\bmu_g).
\]
DEIM treats $(\mathcal{T}(x; \bmu_g))_{k,\ell}$ with $x$ discretized, thus approximating the $\bmu_g$-dependent $\calN \times 1$ vector by a linear combination of $\bmu_g$-independent $\calN \times 1$ vectors. 
This means that, to obtain an affine approximation of $\mathcal{T}(x; \bmu_g)$, we need $Q_a = \sum_{k, \ell =1}^2 Q_{k \ell}$ terms. Each $Q_{k \ell}$, and thus the total $Q_a$, tend to be large when the parameterized operator involves geometrical deformations \cite{Benaceur2018, Chen2012}. 
As a consequence, the online efficiency, which is dependent on $Q_a$, is severely degraded. 
The matrix version of discrete empirical interpolation method (MDEIM) \cite{negri2015efficient} partially alleviate the situation by treating the matrix $\mathcal{K}(\bmu_g)$ as a whole. That is, it finds $Q_a$ parameter-independent matrices ${\mathcal K}_1, \dots, {\mathcal K}_{Q_a}$
\[
{\mathcal K}(\bmu_g) = \sum_{q=1}^{Q_a} {\mathcal K}_q \theta_q(\bmu_g).
\]
{This decomposition} is achieved by expressing the matrix ${\mathcal K}(\bmu_g)$ in vector format via stacking its columns, and then applying DEIM to the resulting ($\bmu_g$-dependent) vector.
}

{However, $Q_a$ is large even with MDEIM \cite{negri2015efficient}. It is therefore critical to test our newly proposed hybrid methods together with the other offline-enhanced greedy algorithms on this particular problem. Toward that end, we test two different cases. The first case includes two parameters: frequency and one RBF control point, and the second one contains all five parameters $\begin{bmatrix}f & \bmu_g \end{bmatrix}$.  {MDEIM \cite{negri2015efficient} is applied to obtain the affine approximation, as is done therein.} In the numerical experiments,} ${\mathbb {TSD\_CG}}$ and its corresponding hybrid ${\mathbb {H\_TSD\_CG}}$ fail {for both cases}. As a result, we report the comparison of {the four algorithms that work, that is} ${\mathbb {CG}}$, ${\mathbb {STS\_CG}}$, ${\mathbb {AE\_CG}}$ and ${\mathbb {H\_AE\_CG}}$. {The detailed setup of both cases are provided in Table \ref{tab:CompDetail}.}

\begin{table}[h!]
  \begin{center}
  \resizebox{\textwidth}{!}{
    \renewcommand{\tabcolsep}{1cm}
    \renewcommand{\arraystretch}{1}
 {
\begin{tabular}{@{}l|p{0.22\textwidth}@{}|p{0.22\textwidth}@{}}
      \toprule    
Variable & Value (two parameters) & Value (five parameters)\\
\midrule
Size of {$\Xi_{\rm train}$}  & $10000$&$10000$\\   
Number of vertices &4567& 4567\\
Number of elements & 8740 & 8740\\
Number of nodes & 4567& 4567\\
Number of finite element dofs & 4567& 4567\\
Interpolation procedure in (\ref{eq:affine})
&  MDEIM\cite{negri2015efficient}  & MDEIM\cite{negri2015efficient} \\
Interpolation points used in (\ref{eq:affine}) & 30 & 250\\
{Number of matrix bases} $Q_a$ in (\ref{eq:affine}) & 3 & 95\\
{Number of RHS bases} $Q_f$ in (\ref{eq:affine}) & 1 &4\\
$(K_{\rm damp}, C_M, M_{\rm sample})$ & ( {5}, 10, 2048)&( {5},10, 2048)\\
Tolerance $\varepsilon_{\rm tol}$  &{various}     & {various} \\
 \bottomrule
    \end{tabular}
}
\renewcommand{\arraystretch}{1}
    \renewcommand{\tabcolsep}{12pt}
}
  \end{center}
  \caption{Experiment setup for the acoustic horn.}
  \label{tab:CompDetail}
\end{table}

\subsubsection{Two parameters case}
We consider frequency and the vertical displacement of the right-most RBF control point in Figure \ref{fig:hornsetup}. The parameter domain is given by $\mathcal{D} = [50, 1000] \times [-0.03, 0.03]$ and $\Xi_{\rm train}$ is generated by operating Latin hypercube sampling in $\mathcal{D}$. Other computational parameters are listed in Table \ref{tab:CompDetail} {second} column.

\begin{figure}[htbp]
\centering
\includegraphics[width=0.49\textwidth]{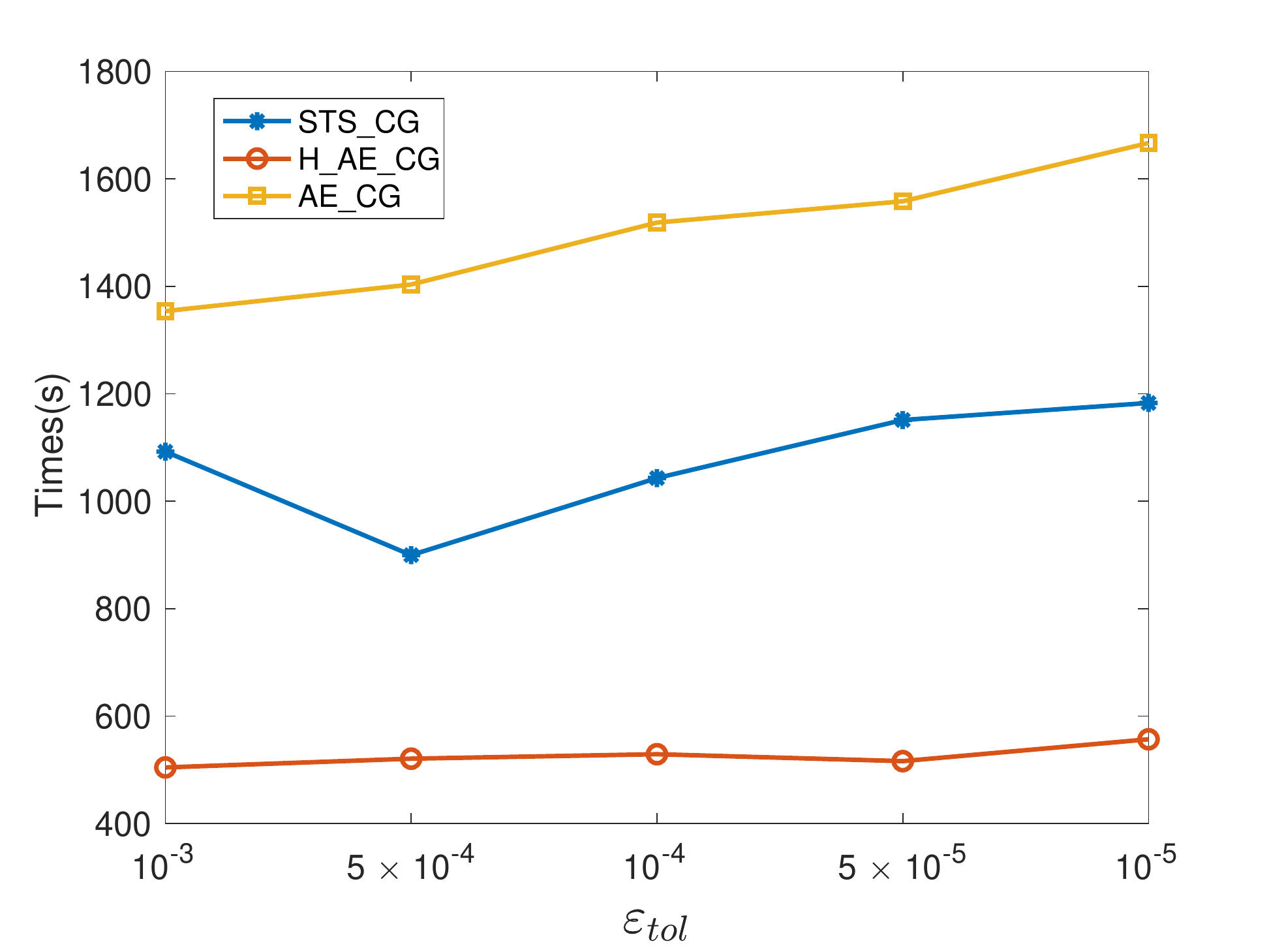}
\includegraphics[width=0.49\textwidth]{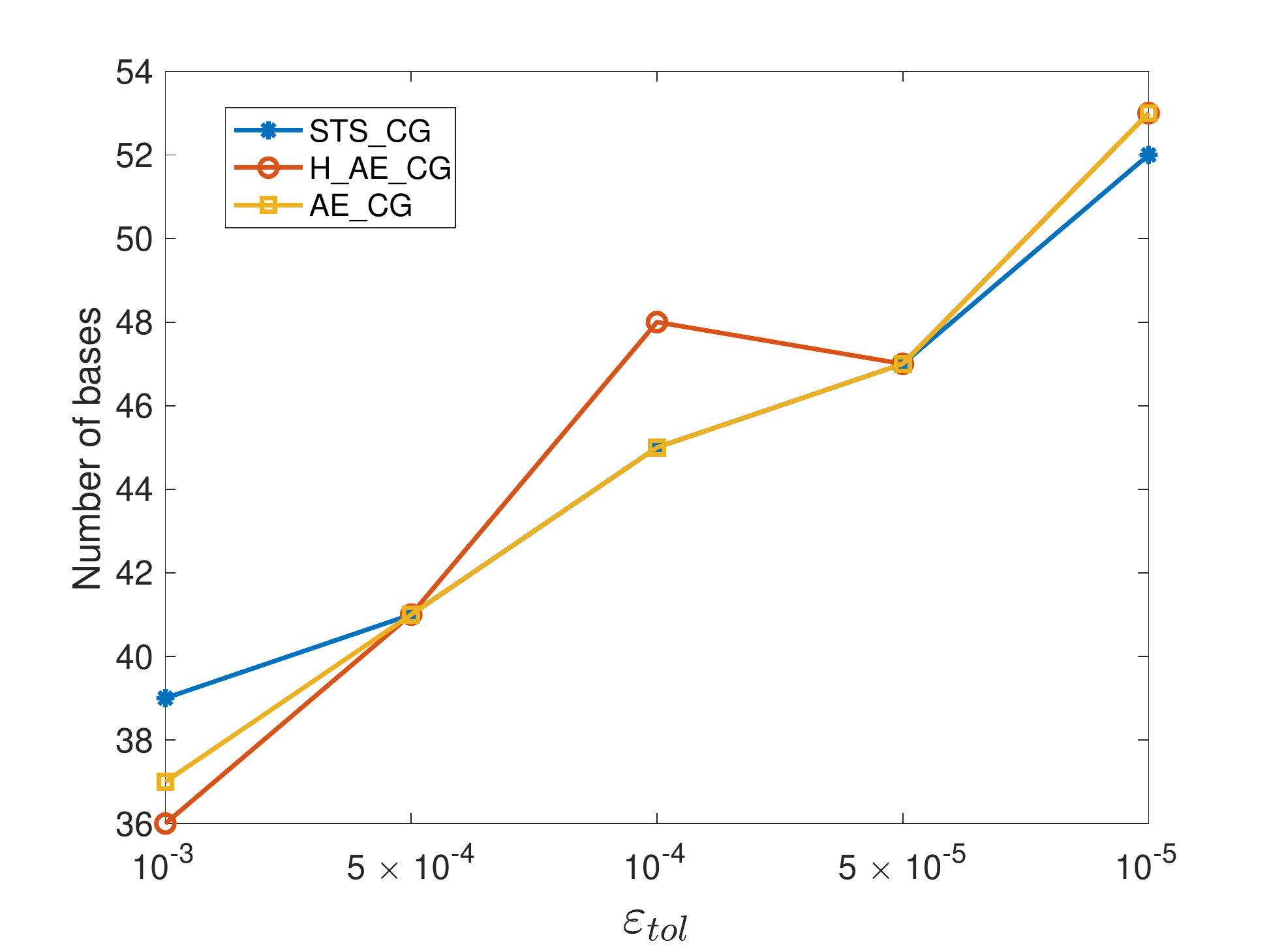}
\caption{Two parameter case for the acoustic horn with different tolerance: computational runtime (left) and number of bases (right).}
\label{fig:time-num}
\end{figure}

Figure \ref{fig:time-num} demonstrates the performance of ${\mathbb {STS\_CG}}$, ${\mathbb {AE\_CG}}$ and ${\mathbb {H\_AE\_CG}}$ for five different tolerances $\varepsilon_{\rm tol} = 10^{-3}, 5 \times 10^{-4}, 10^{-4}, 5\times 10^{-5}, 10^{-5}$. For ${\mathbb {AE\_CG}}$ and ${\mathbb {STS\_CG}}$, smaller tolerance leads to more runtime.  This is much less severe for ${\mathbb {H\_AE\_CG}}$.  This observation suggests that ${\mathbb {H\_AE\_CG}}$ has great potential to handle larger-scale problems. In addition, the sizes of these three RB spaces are very similar, indicating that ${\mathbb {H\_AE\_CG}}$ does not appear to suffer from online efficiency degradation for this example.

\subsubsection{Five parameters case}
We now let all five parameters $\begin{bmatrix}f  & \bmu_{g}\end{bmatrix}$ vary. The parameter domain is given by $\mathcal{D} = [50, 1000] \times \mathcal{D}_{g}$, where $\mathcal{D}_{g} = [-0.03, 0.03]^4$. The full training set $\Xi_{\rm train}$ is also from Latin hypercube sampling in $\mathcal{D}$. The corresponding computational settings are detailed in Table \ref{tab:CompDetail} column 3.

This five-parameter case is highly non-affine and much more complicated than the previous tests. {$95$ matrix bases and $4$ vector bases are required to approximate this non-affine PDE operator, which significantly increases the computational burden for assembling the reduced solver and evaluating the {\it a posteriori error estimate}.} Due to the computation of ${\mathbb {CG}}$ being much more demanding for this case, we test 3 cases: $\varepsilon_{\rm tol} = 10^{-3}, 5\times10^{-4}, 10^{-4}$ to verify the performance of ${\mathbb {STS\_CG}}$, ${\mathbb {AE\_CG}}$ and ${\mathbb {H\_AE\_CG}}$. 
The results are shown in Table \ref{tab:5-para} which demonstrates the efficacy of the new hybrid approach ${\mathbb {H\_AE\_CG}}$.

\begin{table}
\begin{tabularx}{\linewidth}{X X X X X }
    \multicolumn{5}{c}{(a) Runtime(s)} \\
    \toprule
    \multicolumn{2}{c}{Methods/$\varepsilon_{\rm tol}$} & $10^{-3}$ & $5\times10^{-4}$ & $10^{-4}$\\
    \midrule
    \multicolumn{2}{c}{${\mathbb {STS\_CG}}$} & 33849 & 42903 &  78998 \\
    \multicolumn{2}{c}{${\mathbb {AE\_CG}}$} & 33779 & 47473 &  114030 \\
    \multicolumn{2}{c}{${\mathbb {H\_AE\_CG}}$} & 27522 & 37466 &  51003 \\
    \bottomrule
 \end{tabularx}
 
 \bigskip
 
 \begin{tabularx}{\linewidth}{X X X X X }
    \multicolumn{5}{c}{(a) Number of bases} \\
    \toprule
    \multicolumn{2}{c}{Methods/$\varepsilon_{\rm tol}$} & $10^{-3}$ & $5\times10^{-4}$ & $10^{-4}$\\
    \midrule
    \multicolumn{2}{c}{${\mathbb {STS\_CG}}$} & 109 & 130 &  193 \\
    \multicolumn{2}{c}{${\mathbb {AE\_CG}}$} & 105 & 127 &  187 \\
    \multicolumn{2}{c}{${\mathbb {H\_AE\_CG}}$} & 116 & 135 & 199  \\
    \bottomrule
 \end{tabularx}
 \caption{Results for the five-parameter acoustic horn with different tolerances.}
 \label{tab:5-para}
\end{table}

\section{Conclusions}

\label{sec5}
{One of the main computational bottleneck in developing reduced basis methods for parametrized PDEs is the deteriorating offline efficiency for problems with large-scale training set or high-dimensional parameter space. While there are several past attempts to pursue this, these methods still suffer from significant computational cost, and/or lack of robustness resulting from random sampling, and/or sensitivities to the algorithms' tuning parameters.}

In this paper, we review three recent {such} offline enhancement approaches for the reduced basis method {all of} {which} share the  overarching theme of constructing a small-size subset of the full training set and then performing the classical greedy algorithm on it. {By closely integrating them, we obtain two new hybrid approaches which are remarkably faster without impacting the quality of reduced basis space. Through two types of different numerical experiments, we demonstrate a significant cost reduction for the offline portion of RBM}, robustness, {and} ease of setting tuning parameters. {These experiments lead us to believe that the efficiency will be more pronounced when the non-affinity gets more severe, or the dimension of parameter space gets higher.  A detailed follow-up study constitutes our future work.}

%\backmatter

\section*{Acknowledgments}
The authors want to express their sincere gratitude to Prof. Anthony Patera of MIT for suggesting to test the algorithms on the Helmholtz equation. The second author was partially supported by National Science Foundation grant DMS-1719698. This material is based upon work supported by the National Science Foundation under Grant No. DMS-1439786 and by the Simons Foundation Grant No. 50736 while Y.~Chen and J.~Jiang were in residence at the Institute for Computational and Experimental Research in Mathematics in Providence, RI, during the ``Model and dimension reduction in uncertain and dynamic systems" program.

\bibliography{ultra_RB}

\end{document}